\documentclass{elsarticle}

\usepackage{amsmath}
\usepackage{amssymb}
\usepackage{gastex}
\usepackage{graphicx}
\usepackage[usenames]{color}
\usepackage{xspace}
\usepackage{float}

\usepackage{stmaryrd}

\newtheorem{Theorem}{Theorem}
\newtheorem{Lemma}[Theorem]{Lemma}
\newtheorem{Proposition}[Theorem]{Proposition}
\newtheorem{Corollary}[Theorem]{Corollary}
\newdefinition{Definition}[Theorem]{Definition}
\newdefinition{Example}[Theorem]{Example}
\newdefinition{Remark}[Theorem]{Remark}
\newproof{Proof}{Proof}
\newenvironment{QProof}{\begin{Proof}}{\qed\end{Proof}}
\newproof{Proof-app}{Proof of Proposition~\ref{prop:q1_2_Aigner}}

\newcommand{\two}{\mathrm{GF}(2)}
\newcommand{\gffour}{\mathrm{GF}(4)}

\newcommand{\sub}[1]{[ #1 ]}

\newcommand{\sdif}{\mathop{\mathrm{\Delta}}}

\newcommand{\vertexrem}{\setminus}
\newcommand{\setminor}{\setminus}

\newcommand{\dmatroid}{$\Delta$-matroid\xspace}
\newcommand{\dmatroids}{$\Delta$-matroids\xspace}
\newcommand{\dual}{\mathop{\bar{*}}}
\newcommand{\vfsafe}{vf-safe\xspace}

\renewcommand{\emptyset}{\varnothing}

\begin{document}

\begin{frontmatter}

\title{Interlace Polynomials for Multimatroids and Delta-Matroids}

\author[rb]{Robert Brijder\corref{cor}}
\cortext[cor]{Corresponding author}
\ead{robert.brijder@uhasselt.be}

\author[hjh]{Hendrik Jan Hoogeboom}

\address[rb]{Hasselt University and Transnational University of Limburg, Belgium}
\address[hjh]{Leiden Institute of Advanced Computer Science,\\
Leiden University, The Netherlands}


\begin{abstract}
We provide a unified framework in which the interlace polynomial and several related graph polynomials are defined more generally for multimatroids and delta-matroids. Using combinatorial properties of multimatroids rather than graph-theoretical arguments, we find that various known results about these polynomials, including their recursive relations, are both more efficiently and more generally obtained. In addition, we obtain several interrelationships and results for polynomials on multimatroids and delta-matroids that correspond to new interrelationships and results for the corresponding graphs polynomials. As a tool we prove the equivalence of tight 3-matroids and delta-matroids closed under the operations of twist and loop complementation, called vf-safe delta-matroids. This result is of independent interest and related to the equivalence between tight 2-matroids and even delta-matroids observed by Bouchet.
\end{abstract}

\end{frontmatter}

\section{Introduction}
The discovery of the interlace polynomial
by Arratia, Bollob\'as, and Sorkin
\cite{Arratia_InterlaceP_SODA,Arratia2004199} can be traced from 4-regular graphs, via circle graphs, to general graphs. In this paper we generalize this graph polynomial to \dmatroids and, even more generally, to multimatroids. In fact, we introduce generic polynomials for both \dmatroids and multimatroids that allow for a unified presentation of the interlace polynomial and various related graph polynomials. This introduction provides a brief sketch of a number of these known polynomials and their relationships, see, e.g.,
\cite{Ellis/01/SurveyGraphPolsI,Ellis/01/SurveyGraphPolsII} for a survey.

The Martin polynomial \cite{MartinPoly1977}, which is defined
for both directed and undirected graphs, computes the number of
circuit partitions of that graph. More precisely, the number
$a_k$ of $k$-component circuit partitions in a graph $D$ is
equal to the coefficient of $y^k$ of the Martin polynomial of
$D$. In case $D$ is a 2-in, 2-out digraph, an Eulerian circuit $C$
of $D$ defines a circle graph $G$. The interlace polynomial of
$G$ turns out to coincide (up to a trivial transformation) with
the Martin polynomial of (the underlying 2-in, 2-out digraph)
$D$ --- the interlace polynomial is however defined for graphs
in general (i.e., not only for circle graphs)
\cite{Arratia_InterlaceP_SODA,Arratia2004199}. The interlace
polynomial is known to be invariant under the graph operations
of local and edge complementation (where local complementation
is defined here only for looped vertices) and moreover the
interlace polynomial fulfills a recursive reduction relation.
As demonstrated in \cite{Aigner200411,ArratiaBS04}, the
interlace polynomial $q(G)$ of a graph $G$ may be explicitly
defined through the nullity values of the set of subgraphs of
$G$. The interlace polynomial moreover coincides with the
Tutte-Martin polynomial defined for isotropic
systems~\cite{Arratia2004199,DBLP:journals/dm/Bouchet05} in the
case where graph $G$ does not have loops. Restricting to such
simple graphs is an essential loss of generality as, e.g., the
Tutte-Martin polynomial does not explain the invariance of the
interlace polynomial under local complementation.

Related polynomials have been defined such as the different
``interlace polynomial'' $Q(G)$ of \cite{Aigner200411}, which is
invariant under local, loop, and edge complementation. In case $G$
is a circle graph (with possibly loops), then $Q(G)$ coincides with
the Martin polynomial of a 4-regular graph corresponding to $G$.
Moreover, $Q(G)$ coincides with the ``global'' Tutte-Martin
polynomial defined for isotropic
systems~\cite{DBLP:journals/dm/Bouchet05} in the case where graph
$G$ is simple. As $Q(G)$ is invariant under loop complementation,
this is not an essential loss of generality. Also, the bracket
polynomial for graphs \cite{Traldi/Bracket1/09} has recursive
relations similar to the recursive relations of the interlace
polynomial.

The notion of delta-matroid (or \dmatroid), introduced by Bouchet
\cite{bouchet1987}, is a generalization of the notion of matroid. In
addition, binary \dmatroids may be viewed as a generalization of graphs.
Moreover, as pointed out in \cite{Geelen97}, the graph operations of
local and edge complementation have particularly simple
interpretations in terms of \dmatroids. In
\cite{BH/PivotLoopCompl/09} a suitable generalization of the notion
of loop complementation for graphs has been defined for a subclass of the \dmatroids called vf-safe \dmatroids.
Moreover, in \cite{BH/NullityLoopcGraphsDM/10}, the notion of
nullity for the adjacency matrix of a graph (or skew-symmetric matrix in general) is shown to correspond to a natural distance
measure within \dmatroids.

Multimatroids were introduced by Bouchet \cite{DBLP:journals/siamdm/Bouchet97} as a common generalization of both \dmatroids and isotropic systems. It turns out that \dmatroids precisely correspond to multimatroids with two ``directions'', called $2$-matroids. A particularly interesting class of multimatroids are the tight multimatroids, defined in \cite{DBLP:journals/ejc/Bouchet01}. The so-called even \dmatroids precisely correspond to the tight $2$-matroids, while the class of isotropic systems corresponds to a subclass of tight $3$-matroids (i.e., tight multimatroids with three directions). We show in Section~\ref{sec:equivalence_vf-safe} that the class of vf-safe \dmatroids (those \dmatroids that ``allow'' loop complementation in conjunction with twist) precisely corresponds to the class of tight $3$-matroids.
As a consequence, vf-safe \dmatroids strictly generalize isotropic systems.
The equivalence of vf-safe \dmatroids and tight $3$-matroids is of independent interest and confirms that the vf-safe \dmatroids form a natural family even beyond their convenient closure properties.

In Section~\ref{sec:mm_pols} we define the (weighted) transition polynomial $Q(Z)$ for multimatroids $Z$ in a natural way as a ``nullity generating'' polynomial, and moreover prove a number of evaluations of $Q(Z)$ in case $Z$ is a tight $k$-matroid (i.e., with $k$ directions).
Then, in Section~\ref{sec:def_sets_interlacep}, we define the (weighted) transition polynomial $Q(M)$ for set systems $M$ using the generalizations of the notions of nullity, and of local, loop, and edge complementation mentioned above. In Section~\ref{sec:ss_polynomials_recursion} we show that the transition polynomial $Q(Z)$ for tight $3$-matroids $Z$ corresponds to the transition polynomial $Q(M)$ for vf-safe \dmatroids $M$.
It turns out that the various graph polynomials recalled above are all special cases of the transition polynomials for (vf-safe) \dmatroids and multimatroids. In this way, these polynomials provide unified views in which invariance properties, recursive relations, and evaluations of these polynomials may be investigated.

The proofs of these results use combinatorial properties of multimatroids, and are not only more general and unified, but also much more concise compared to the proofs which rely on graph-theoretical arguments. In particular, the evaluation on $-1$ of the interlace polynomial for graphs (and also of the Tutte polynomial for binary matroids) is strikingly simple when viewed as a special case of an evaluation of the transition polynomial for tight multimatroids.
Moreover, due to the highly symmetrical nature of multimatroids, we naturally obtain new interrelationships and evaluations for the corresponding polynomials on graphs (Section~\ref{sec:graph_poly}).

Since the minor-closed class of \vfsafe \dmatroids strictly contains the class of
quaternary matroids (this is recalled in
Section~\ref{sec:notions_DM}), the main results
of this paper are not likely obtainable using isotropic systems
which fundamentally deal with binary matroids.

Inspired by the connection between the Martin and the Tutte
polynomial in \cite{MartinPoly1977} (and its relation to
isotropic systems described in
\cite{TutteMartinOrientingVectors/Bouchet91}) we find that the
Tutte polynomial $t_M(x,y)$ for matroids $M$ can, for the case
$x=y$, be seen as a special case of the transition
polynomial for \dmatroids. In addition, the recursive relations
of the transition polynomial and $t_M(y,y)$ coincide
when restricting to matroids $M$. The obtained evaluations of
the transition polynomial are then straightforwardly
carried over to $t_M(y,y)$.

\setlength{\tabcolsep}{15pt}
\begin{table}
\begin{center}
\begin{tabular}{l|l}
& Special cases\\[0.1cm]
\hline
$Q_{[1,1,0]}(M)(y)$ & Martin polynomial for $2$-in, $2$-out digraphs \cite{MartinPoly1977}, \\
& Single-variable interlace polynomial \cite{Arratia2004199}, \\
& Tutte-Martin polynomial for isotropic systems \cite{TutteMartinOrientingVectors/Bouchet91}, \\
& Tutte polynomial for matroids restricted to the diagonal \\[0.15cm]
$Q_{[1,x,0]}(M)(y)$ & Two-variable interlace polynomial \cite{ArratiaBS04}\\[0.15cm]
$Q_{[a,b,0]}(M)(y)$ & Tutte polynomial for matroids restricted to part of the plane \\[0.15cm]
$Q_{[0,1,x]}(M)(y)$ & Bracket polynomial for graphs \cite{Traldi/Bracket1/09}\\[0.15cm]
$Q_{[1,1,1]}(M)(y)$ & Martin polynomial for $4$-regular graphs \cite{MartinPoly1977},\\
& ``Global'' Tutte-Martin polynomial for isotropic systems \cite{TutteMartinOrientingVectors/Bouchet91}, \\
& Polynomial $Q(G,x)$ for simple graphs $G$ as defined in \cite{Aigner200411} \\
\end{tabular}
\end{center}
\caption{Specializations of $Q(M)$  for (vf-safe) \dmatroids $M$.}
\label{table:summary}
\end{table}
Table~\ref{table:summary} summarizes the relations of the above recalled polynomials from the literature with the transition polynomial $Q(M)$ for (vf-safe) \dmatroids $M$.

\section{Multimatroids}

We assume the reader is familiar with the basic notions concerning matroids, which can be found, e.g., in \cite{Welsh/MatroidBook,Oxley/MatroidBook-2nd}.

We take the terminology of multimatroids as developed by Bouchet \cite{DBLP:journals/siamdm/Bouchet97,DBLP:journals/combinatorics/Bouchet98,DBLP:journals/ejc/Bouchet01}.
A \emph{carrier} is a tuple $(U,\Omega)$ where $\Omega$ is a partition of a finite set $U$, called the \emph{ground set}. Every $\omega \in \Omega$ is called a \emph{skew class}, and a $p \subseteq \omega$ with $|p|=2$ is called a \emph{skew pair} of $\omega$. A \emph{transversal} (\emph{subtransversal}, resp.) $T$ of $\Omega$ is a subset of $U$ such that $|T \cap \omega| = 1$ ($|T \cap \omega| \leq 1$, resp.) for all $\omega \in \Omega$. We denote the set of transversals of $\Omega$ by $\mathcal{T}(\Omega)$, and the set of subtransversals of $\Omega$ by $\mathcal{S}(\Omega)$. The power set of a set $X$ is denoted by $2^X$.

We recall now the notion of multimatroid. Like matroids, multimatroids can be defined in terms of rank, circuits, independent sets, etc. We define multimatroids here in terms of independent sets.
\begin{Definition} [\cite{DBLP:journals/siamdm/Bouchet97}]\label{def:multimatroid}
A multimatroid $Z$ (described by its independent sets) is a triple $(U,\Omega,\mathcal{I})$, where $(U,\Omega)$ is a carrier and $\mathcal{I} \subseteq \mathcal{S}(\Omega)$ such that:
\begin{enumerate}
\item for each $T \in \mathcal{T}(\Omega)$, $(T,\mathcal{I}\cap 2^T)$ is a matroid (described by its independent sets) and
\item for any $I \in \mathcal{I}$ and any skew pair $p = \{x,y\}$ of some $\omega \in \Omega$ with $\omega \cap I = \emptyset$, $I \cup \{x\} \in \mathcal{I}$ or $I \cup \{y\} \in \mathcal{I}$.
\end{enumerate}
\end{Definition}
Terminology concerning $\Omega$ carries over to $Z$: hence we may, e.g., speak of a transversal of $Z$. We often simply write $U$ and $\Omega$ to denote the ground set and the partition of the multimatroid $Z$ under consideration, respectively. Each $I \in \mathcal{I}$ in the definition of multimatroid $Z$ is called an \emph{independent set} of $Z$. The family $\mathcal{I}$ of independent sets of $Z$ is denoted by $\mathcal{I}_Z$. The family $\max(\mathcal{I}_Z)$ of sets of $\mathcal{I}_Z$ that are maximal with respect to inclusion is denoted by $\mathcal{B}_Z$. The elements of $\mathcal{B}_Z$ are called the \emph{bases} of $Z$. Note that the bases uniquely determine $Z$ (given its carrier). For any $X \subseteq U$, the \emph{restriction} of $Z$ to $X$, denoted by $Z[X]$, is the multimatroid $(X,\Omega',\mathcal{I}\cap 2^{X})$ with $\Omega' = \{ \omega \cap X \mid \omega \cap X \neq \emptyset, \omega \in \Omega\}$.
If $X$ is a subtransversal, then we identify $Z[X]$ with the matroid $(X,\mathcal{I}\cap 2^{X})$ since $\Omega' = \{\{u\} \mid u \in X\}$ captures no additional information. The \emph{rank} of $S \in \mathcal{S}(\Omega)$ in $Z$, denoted by $r_Z(S)$, is the rank $r(Z[S])$ of the matroid $Z[S]$. The \emph{nullity} of $S \in \mathcal{S}(\Omega)$ in $Z$, denoted by $n_Z(S)$, is $n_Z(S) = |S| - r_Z(S)$, i.e., the nullity $n(Z[S])$ of the matroid $Z[S]$. The second condition of Definition~\ref{def:multimatroid} can equivalently be formulated as follows: For all $S \in \mathcal{S}(\Omega)$ and any skew pair $p = \{x,y\}$ of $\omega \in \Omega$ with $w \cap S = \emptyset$, $n_Z(S \cup \{x\}) - n_Z(S) + n_Z(S \cup \{y\}) - n_Z(S) \leq 1$ (see the definition of multimatroid in terms of rank in \cite{DBLP:journals/siamdm/Bouchet97}).

The \emph{minor} of $Z$ induced by $X$, denoted by $Z|X$, is the multimatroid $(U',\Omega',\mathcal{I}')$, where $\Omega' = \{\omega \in \Omega \mid \omega \cap X = \emptyset\}$, $U' = \cup_{\omega \in \Omega'} \omega$, and $r_{Z|X}(S) = r_Z(S\cup X) - r_Z(X)$ for all $S \in \mathcal{S}(\Omega')$. By a standard property of matroid contraction (see, e.g., \cite[Proposition~3.1.6]{Oxley/MatroidBook-2nd}), for each
$T \in \mathcal{T}(\Omega')$, the matroid $Z|X[T]$ is equal to $Z[T \cup X] \slash X$, where $\slash$ denotes matroid contraction. In case $X = \{u\}$ is a singleton, we also write $Z|u$ to denote $Z|\{u\}$. An element $u \in U$ is called \emph{singular} in $Z$ if $n_Z(\{u\}) = 1$. Thus if $u \in U$ is nonsingular in $Z$, then $n_Z(\{u\}) = 0$. Note that, by the second condition of Definition~\ref{def:multimatroid}, each skew class contains at most one singular element. A skew class that contains a singular element is called \emph{singular}. We recall the following result of \cite{DBLP:journals/ejc/Bouchet01}.

\begin{Proposition}[Proposition~5.5 of \cite{DBLP:journals/ejc/Bouchet01}] \label{prop:minor_singular}
If a skew class $\omega$ of a multimatroid $Z$ is singular, then $Z|u = Z[U\setminus \omega]$ for all $u \in \omega$.
\end{Proposition}

For any positive integer $q$, a \emph{$k$-matroid} is a multimatroid with carrier $(U,\Omega)$ such that $|\omega|=k$ for all $\omega \in \Omega$. Note that a (regular) matroid corresponds to a $1$-matroid.

A multimatroid $Z$ is called \emph{nondegenerate} if $|\omega| > 1$ for all $\omega \in \Omega$. It is shown in \cite{DBLP:journals/siamdm/Bouchet97}, that if $Z$ is nondegenerate, then $\mathcal{B}_Z = \mathcal{I}_Z \cap \mathcal{T}(\Omega)$.

A multimatroid $Z$ is called \emph{tight} if both $Z$ is nondegenerate and for all $S \in \mathcal{S}(\Omega)$ with $|S| = |\Omega|-1$, there is a $x \in \omega$ such that $n_Z(S\cup \{x\}) = n_Z(S)+1$, where $\omega \in \Omega$ is the unique skew class such that $S \cap \omega = \emptyset$. Consequently, by the second condition in Definition~\ref{def:multimatroid}, the values $n_Z(S\cup \{y\})$ for all $y \in \omega$ are equal, except for exactly one $x \in \omega$ (for which the value $n_Z(S\cup \{x\})$ is one larger than the others). It is shown in \cite[Proposition~4.1]{DBLP:journals/ejc/Bouchet01} that tightness is preserved under minors. Moreover, the notion of tightness is characterized in terms of bases.
\begin{Proposition}[Theorem~4.2c of \cite{DBLP:journals/ejc/Bouchet01}] \label{prop:base-tight}
Let $Z$ be a nondegenerate multimatroid. Then $Z$ is tight iff for every basis $X$ and every skew class $\omega$ of $Z$, exactly one of the transversals $(X \setminus \omega) \cup \{u\}$, $u\in \omega$, is not a basis of $Z$.
\end{Proposition}

As an example, we briefly recall from \cite{DBLP:journals/siamdm/Bouchet97} that the Eulerian circuits of a connected $4$-regular graph form the bases of a tight $3$-matroid.
\begin{Example}
Given a connected $4$-regular graph $G$, any Eulerian circuit $C$ of $G$ visits every vertex of $G$ exactly twice. The Eulerian circuit $C$ defines a \emph{transition} at each vertex $v$, i.e., a set of two (disjoint) unordered pairs of half-edges incident to $v$, that coincides with the trajectory of $C$ at $v$. Note that each vertex has three possible transitions. Let $(U,\Omega)$ be a carrier, where each skew class of $\Omega$ consists of the three possible transitions at a vertex of $G$. The transversals of $(U,\Omega)$ are in one-to-one correspondence with the \emph{circuit partitions} of $G$, where a circuit partition is a set of mutually edge-disjoint circuits (circuits are considered unoriented) of $G$ covering each edge of $G$. It turns out that the transversals of $(U,\Omega)$ that correspond to the Eulerian circuits of $G$ are the bases of a $3$-matroid $Z$, called the \emph{Eulerian multimatroid} of $G$ \cite{DBLP:journals/siamdm/Bouchet97}. It is easily seen by Proposition~\ref{prop:base-tight} that $Z$ is tight: for every Eulerian circuit $C$ of $G$ and every vertex $v$ of $G$, it is possible to split $C$ in two by changing the transition of $C$ at $v$ in a suitable way.
\end{Example}

\section{Polynomials for Multimatroids} \label{sec:mm_pols}

It is often worthwhile to study weighted/multivariate variants of enumerating polynomials, see, e.g., the weighted variants of the Tutte polynomial studied in  \cite{BollobasRiordan:TutteColoured,EllisMonaghanTraldi2006,Sokal-multivariateTutte,ZaslavskyStrongTutte}.
We now define a weighted polynomial enumerating the nullity values of transversals of a multimatroid. It is called the transition polynomial, motivated by the eponymous polynomial for $4$-regular graphs studied by Jaeger \cite{Jaeger1990Transition}.

We obtain in this section a recursive relation for our transition polynomial as well as evaluations of the polynomial for tight multimatroids, following the lead of \cite{Jaeger1990Transition}.

\begin{Definition}
Let $Z$ be a multimatroid with carrier $(U,\Omega)$. We define the \emph{(weighted) transition polynomial} of $Z$ as
$$
Q(Z)(\vec{x}, y) = \sum_{T \in \mathcal{T}(\Omega)} x_T y^{n_Z(T)},
$$
where $\vec{x}$ is a vector indexed by $U$ with entries $x_{u}$ for all $u \in U$, and $x_T = \prod_{u \in T} x_{u}$.
\end{Definition}

We now provide recursive relations for $Q(Z)$.
\begin{Theorem} \label{thm:mm_recursive}
Let $Z$ be a multimatroid and let $\omega \in \Omega$.

If $\omega$ is nonsingular in $Z$, then
$$
Q(Z)(\vec{x}, y) = \sum_{v \in \omega} x_v Q({Z|v})(\vec{x'}, y),
$$
and if $\omega$ is singular in $Z$, then for all $w \in \omega$
$$
Q(Z)(\vec{x}, y) = \left(x_u y + \sum_{v \in \omega \setminus \{u\}} x_v \right) Q({Z|w})(\vec{x'}, y),
$$
where $u \in \omega$ is singular in $Z$, and $\vec{x'}$ is obtained from $\vec{x}$ by removing the entries indexed by $\omega$.
\end{Theorem}
\begin{QProof}
Let $\Omega' = \Omega \setminus \{\omega\}$.

We have
\begin{eqnarray*}
Q(Z)(\vec{x}, y) & = & \sum_{u \in \omega} \; \sum_{\substack{T \in \mathcal{T}(\Omega)\\ u \in T}} x_T y^{n_Z(T)} = \sum_{u \in \omega} x_{u} \sum_{\substack{T \in \mathcal{T}(\Omega)\\ u \in T}} x_{T\setminus \{u\}} y^{n_Z(T)}\\
& = & \sum_{u \in \omega} x_{u} \sum_{T \in \mathcal{T}(\Omega')} x_{T} y^{n_Z(T\cup\{u\})}.
\end{eqnarray*}
Let $u \in \omega$. We have, for all $S \in \mathcal{S}(\Omega')$, $r_{Z|u}(S) = r_Z(S \cup \{u\}) - r_Z(\{u\})$ or, equivalently, $n_{Z|u}(S) = n_Z(S \cup \{u\}) - n_Z(\{u\})$.

Assume first that $u$ is nonsingular in $Z$. Then $n_{Z|u}(S) = n_Z(S \cup \{u\})$ for all $S \in \mathcal{S}(\Omega')$, and thus $\sum_{T \in \mathcal{T}(\Omega')} x_{T} y^{n_Z(T\cup\{u\})} = Q({Z|u})(\vec{x'}, y)$.

Assume now that $u$ is singular in $Z$. Then $n_{Z|u}(S) = n_Z(S \cup \{u\}) - 1$ for all $S \in \mathcal{S}(\Omega')$, and thus $\sum_{T \in \mathcal{T}(\Omega')} x_{T} y^{n_Z(T\cup\{u\})} = y Q({Z|u})(\vec{x'}, y)$. The result follows by Proposition~\ref{prop:minor_singular}.
\end{QProof}

Note that the recursive relations of Theorem~\ref{thm:mm_recursive}, along with the fact that $Q(Z)(\vec{x}, y) = 1$ when $Z$ is the unique multimatroid with empty ground set, \emph{characterize} the transition polynomial $Q(Z)(\vec{x}, y)$.

The following result is inspired by \cite[Proposition~10]{Jaeger1990Transition} from the context of $4$-regular graphs.
\begin{Theorem} \label{thm:mm_Q_const}
Let $Z$ be a tight multimatroid, let $\omega \in \Omega$, and let $\vec{x'}$ be a vector indexed by $U$ with $x'_v = x'_w$ for all $v,w \in \omega$ and $x'_v = 0$ for all $v \in U\setminus \omega$. Then $Q(Z)(\vec{x}, 1-k) = Q(Z)(\vec{x}+\vec{x'}, 1-k)$ with $k = |\omega|$.
\end{Theorem}
\begin{QProof}
Let $c = x'_v$ for all $v \in \omega$. We prove by induction.

Assume that $|\Omega|=1$. Since $Z$ is tight, there is a $u \in \omega$ singular in $Z$. By Theorem~\ref{thm:mm_recursive}, $Q(Z)(\vec{x}+\vec{x'},1-k) = (x_u+c)(1-k) + \sum_{v \in \omega \setminus \{u\}} (x_v + c) = c(1-k+k-1) + x_u(1-k) + \sum_{v \in \omega \setminus \{u\}} x_v = Q(Z)(\vec{x},1-k)$.

Assume that $|\Omega|>1$. Let $\omega' \in \Omega$ with $\omega' \neq \omega$. Note that $\vec{x}$ and $\vec{x}+\vec{x'}$ coincide for the elements of $\omega'$. The induction hypothesis ensures now (recall that tightness is preserved under minors) that the right-hand sides of the equalities in Theorem~\ref{thm:mm_recursive} coincide for $Q(Z)(\vec{x},1-k)$ and $Q(Z)(\vec{x}+\vec{x'},1-k)$.
\end{QProof}

The polynomial $Q(Z)(\vec{x}, y)$, where $x_{u} = 1$ for all $u \in U$, is denoted by $Q_1(Z)(y) = \sum_{T \in \mathcal{T}(\Omega)} y^{n_Z(T)}$. Note that $Q_1(Z)(1) = |\mathcal{T}(\Omega)| = \prod_{\omega \in \Omega} |\omega|$ and $Q_1(Z)(0) = |\{ T\in\mathcal{T}(\Omega) \mid n_Z(T)=0 \}| = |\mathcal{B}_Z|$.

We now obtain the following two corollaries to Theorem~\ref{thm:mm_Q_const}.
\begin{Theorem} \label{thm:mm_Q_1-k}
Let $Z$ be a tight multimatroid with $U \not= \emptyset$. If $\omega$ is a skew class of $Z$ with $|\omega| = k$, then $Q_1(Z)(1-k) = 0$.
\end{Theorem}
\begin{QProof}
By Theorem~\ref{thm:mm_Q_const}, $Q_1(Z)(1-k) = Q(Z)(\vec{x}, 1-k)$ with $\vec{x}$ the zero vector indexed by $U$. Since $U \not= \emptyset$, $Q(Z)(\vec{x}, 1-k) = 0$.
\end{QProof}

The following result is inspired by \cite[Proposition~11]{Jaeger1990Transition} from the context of $4$-regular graphs.
\begin{Theorem} \label{thm:mmpol_res_1-k}
Let $Z$ be a tight $k$-matroid for some $k > 1$ and let $T \in \mathcal{T}(\Omega)$. Then $Q_1(Z[U\setminus T])(1-k) = (-1)^{|\Omega|}(1-k)^{n_Z(T)}$.
\end{Theorem}
\begin{QProof}
Note that $Q_1(Z[U\setminus T])(y)$ is equal to the polynomial $Q(Z)(\vec{x}, y)$, where $x_{u} = 1$ if $u \in U\setminus T$ and $x_{u} = 0$ if $u \in T$. Moreover, $Q(Z)(\vec{x'}, 1-k) = (-1)^{|\Omega|}(1-k)^{n_Z(T)}$ where $x'_{u} = 0$ if $u \in U\setminus T$ and $x'_{u} = -1$ if $u \in T$. Since $\vec{x'}$ is obtained from $\vec{x}$ by subtracting $1$ for each component, the result follows from Theorem~\ref{thm:mm_Q_const} (by iteration on every skew class).
\end{QProof}

It is interesting to remark that, by Theorem~\ref{thm:mmpol_res_1-k}, $Q_1(Z')(1-k)$ for a particular $(k-1)$-matroid $Z'$ depends on the nullity of a set $T$ in a larger tight $k$-matroid $Z$, where $T$ is \emph{disjoint} from the ground set of $Z'$. We show in Theorem~\ref{thm:extend_tight_mm} that $Z$ is essentially unique given $Z'$.

\section{Pivot and Loop Complementation for Set Systems and Delta-matroids}
\label{sec:notions_DM}
In this section we recall the algebra of set systems generated
by the operations of pivot and loop complementation
\cite{BH/PivotLoopCompl/09}. Additionally we recall from
\cite{BH/NullityLoopcGraphsDM/10} the notion of distance in set
systems and its main properties with respect to \dmatroids.

A \emph{set system} (over $V$) is a tuple $M = (V,D)$ with $V$ a
finite set called the \emph{ground set} and $D \subseteq 2^V$ a
family of subsets of $V$. Let $X \subseteq V$. We define $M \sub{X}
= (X,D')$ where $D' = \{Y \in D \mid Y \subseteq X\}$, and define $M
\setminor X = M \sub{V \setminus X}$. In case $X = \{u\}$ is a
singleton, we also write $M \vertexrem u$ to denote $M \vertexrem
\{u\}$. Set system $M$ is called \emph{nonempty} if $D \not=
\emptyset$. We write simply $Y \in M$ to denote $Y \in D$. A set
system $M$ with $\emptyset \in M$ is called \emph{normal} (note that
every normal set system is nonempty).
In particular the set system
$(\,\emptyset, \{\emptyset\}\,)$ is normal. A set system $M$ is
called \emph{equicardinal} if for all $X_1, X_2 \in M$, $|X_1| =
|X_2|$. For convenience we will often simply denote the ground set
of the set system under consideration by $V$.

Let $M = (V,D)$ be a set system. We define, for $X \subseteq
V$, \emph{pivot} (also called \emph{twist} in the literature) of $M$ on $X$, denoted by $M * X$, as $(V,D *
X)$, where $D * X = \{Y \sdif X \mid Y \in D\}$ and $\sdif$ denotes symmetric difference. In case $X =
\{u\}$ is a singleton, we also write simply $M * u$. Moreover,
we define, for $u \in V$, \emph{loop complementation} of $M$ on
$u$ (the motivation for this name is from graphs, see
Section~\ref{sec:graphs}), denoted by $M + u$, as
$(V,D')$, where $D' = D \sdif \{X \cup \{u\} \mid X \in D, u
\not\in X\}$. We assume left associativity of set system
operations. Therefore, e.g., $M+u \vertexrem u *v$ denotes
$((M+u)\vertexrem u)*v$.

It has been shown in \cite{BH/PivotLoopCompl/09} that pivot ${}*u$ and
loop complementation ${}+u$ on a common element $u \in V$ are
involutions
(i.e., of order $2$) that generate a group $F_u$ isomorphic to
$S_3$, the group of permutations on $3$ elements. In
particular, we have ${}+u*u+u = {}*u+u*u$, which is the third
involution (in addition to pivot and loop complementation), and
is called \emph{dual pivot}, denoted by $\dual $. Explicitly, if $M = (V,D)$ and $u \in V$, then $M \dual u = (V,D')$, where $D' = D \sdif \{X \setminus \{u\} \mid X \in D, u
\in X\}$. The
elements of $F_u$ are called \emph{vertex flips}.\footnote{The
notion of vertex flip as defined in this paper corresponds to
the notion of invertible vertex flip in
\cite{BH/PivotLoopCompl/09} --- as we consider only invertible
vertex flips in this paper for notational convenience we omit
the adjective ``invertible''.} We have, e.g., ${}+u*u = {}\dual
u+u = {}*u \dual  u$ and ${}*u+u = {}+u \dual u = {}\dual  u
*u$ for $u \in V$ (which are the two vertex flips in $F_u$ of
order $3$).

While on a single element the vertex flips behave as the group
$S_3$, they commute when applied on different elements.
Hence, e.g., $M*u+v= M+v*u$ and
$M\dual u+v= M+v\dual u$ when $u \not= v$. Also, $M+u+v= M+v+u$
and thus we (may) write, for $X = \{u_1,u_2,\ldots,u_n\} \subseteq
V$, $M+X$ to denote $M + u_1 \ldots + u_n$ (as the result is
independent of the order in which the operations ${}+u_i$ are
applied). Similarly, we define $M\dual X$ for $X \subseteq V$. We
will often use the above equalities without explicit mention.

One may explicitly define the sets in $M*V$, $M+V$, and $M\dual
V$ as follows: $X \in M*V$ iff $V-X \in M$, and $X \in M+V$ iff
$|\{ Z\in M \mid Z\subseteq X \}|$ is odd. Dually, $X \in
M\dual V$ iff $|\{ Z\in M \mid X\subseteq Z \}|$ is odd. In
particular, $\emptyset\in M\dual V$ iff the number of sets in
$M$ is odd.

\begin{figure}
\unitlength 0.50mm
\centerline{%
{\begin{picture}(250,200)\unitlength 0.50mm
\node[Nframe=n,Nw=55,Nh=60](EEN)(25,75){\unitlength 0.50mm
\begin{picture}(50,60)(5,0)
\gasset{AHnb=0}
\gasset{Nw=10,Nh=10,Nmr=5,Nframe=y,loopdiam=8}
\node(p)(30,95){$p$}
\node(q)(15,80){$q$}
\node(r)(45,80){$r$}
\drawedge(p,q){}
\drawedge(p,r){}
\drawedge(q,r){}
\drawloop[loopangle=90](p){}
\drawloop[loopangle=-20](r){}
\gasset{Nw=10,Nh=10,Nmr=5,Nframe=n}
\node(011)(35,30){\small$q$}
\node[Nframe=y](010)(35,00){$\emptyset$}
\node[Nframe=y](001)(10,40){\small$qr$}
\node[Nframe=y](000)(10,10){\small$r$}
\node[Nframe=y](111)(50,45){\small$pq$}
\node[Nframe=y](110)(50,15){\small$p$}
\node(101)(25,55){\small$pqr$}
\node(100)(25,25){\small$pr$}
\drawedge(000,001){}\drawedge(010,011){}
\drawedge(000,010){}\drawedge(001,011){}
\drawedge(100,101){}\drawedge(110,111){}
\drawedge(100,110){}\drawedge(101,111){}
\drawedge(100,000){}\drawedge(101,001){}
\drawedge(110,010){}\drawedge(111,011){}
\end{picture}}
\node[Nframe=n,Nw=50,Nh=60](TWEE)(90,120){\unitlength 0.50mm
\begin{picture}(50,60)(5,0)
\gasset{AHnb=0}
\gasset{Nw=10,Nh=10,Nmr=5,Nframe=y,loopdiam=8}
\node(p)(30,95){$p$}
\node(q)(15,80){$q$}
\node(r)(45,80){$r$}
\drawedge(p,q){}
\drawedge(p,r){}
\drawedge(q,r){}
\drawloop[loopangle=200](q){}
\gasset{Nw=10,Nh=10,Nmr=5,Nframe=n}
\node[Nframe=y](011)(35,30){\small$q$}
\node[Nframe=y](010)(35,00){$\emptyset$}
\node[Nframe=y](001)(10,40){\small$qr$}
\node(000)(10,10){\small$r$}
\node[Nframe=y](111)(50,45){\small$pq$}
\node(110)(50,15){\small$p$}
\node[Nframe=y](101)(25,55){\small$pqr$}
\node[Nframe=y](100)(25,25){\small$pr$}
\drawedge(000,001){}\drawedge(010,011){}
\drawedge(000,010){}\drawedge(001,011){}
\drawedge(100,101){}\drawedge(110,111){}
\drawedge(100,110){}\drawedge(101,111){}
\drawedge(100,000){}\drawedge(101,001){}
\drawedge(110,010){}\drawedge(111,011){}
\end{picture}}
\node[Nframe=n,Nw=50,Nh=60](DRIE)(155,120){\unitlength 0.50mm
\begin{picture}(50,60)(5,0)
\gasset{AHnb=0}
\gasset{Nw=10,Nh=10,Nmr=5,Nframe=y,loopdiam=8}
\node(p)(30,95){$p$}
\node(q)(15,80){$q$}
\node(r)(45,80){$r$}
\drawedge(p,q){}
\drawedge(q,r){}
\drawloop[loopangle=90](p){}
\drawloop[loopangle=200](q){}
\drawloop[loopangle=-20](r){}
\gasset{Nw=10,Nh=10,Nmr=5,Nframe=n}
\node[Nframe=y](011)(35,30){\small$q$}
\node[Nframe=y](010)(35,00){$\emptyset$}
\node(001)(10,40){\small$qr$}
\node[Nframe=y](000)(10,10){\small$r$}
\node(111)(50,45){\small$pq$}
\node[Nframe=y](110)(50,15){\small$p$}
\node[Nframe=y](101)(25,55){\small$pqr$}
\node[Nframe=y](100)(25,25){\small$pr$}
\drawedge(000,001){}\drawedge(010,011){}
\drawedge(000,010){}\drawedge(001,011){}
\drawedge(100,101){}\drawedge(110,111){}
\drawedge(100,110){}\drawedge(101,111){}
\drawedge(100,000){}\drawedge(101,001){}
\drawedge(110,010){}\drawedge(111,011){}
\end{picture}}
\node[Nframe=n,Nw=55,Nh=60](VIER)(220,75){\unitlength 0.50mm
\begin{picture}(50,60)(5,0)
\gasset{AHnb=0}
\gasset{Nw=10,Nh=10,Nmr=5,Nframe=y,loopdiam=8}
\node(p)(30,95){$p$}
\node(q)(15,80){$q$}
\node(r)(45,80){$r$}
\drawedge(p,q){}
\drawedge(q,r){}
\gasset{Nw=10,Nh=10,Nmr=5,Nframe=n}
\node(011)(35,30){\small$q$}
\node[Nframe=y](010)(35,00){$\emptyset$}
\node[Nframe=y](001)(10,40){\small$qr$}
\node(000)(10,10){\small$r$}
\node[Nframe=y](111)(50,45){\small$pq$}
\node(110)(50,15){\small$p$}
\node(101)(25,55){\small$pqr$}
\node(100)(25,25){\small$pr$}
\drawedge(000,001){}\drawedge(010,011){}
\drawedge(000,010){}\drawedge(001,011){}
\drawedge(100,101){}\drawedge(110,111){}
\drawedge(100,110){}\drawedge(101,111){}
\drawedge(100,000){}\drawedge(101,001){}
\drawedge(110,010){}\drawedge(111,011){}
\end{picture}}
%
%
\node[Nframe=n,Nw=50,Nh=60](VIJF)(155,30){\unitlength 0.50mm
\begin{picture}(50,60)(5,0)
\gasset{AHnb=0}
\gasset{Nw=10,Nh=10,Nmr=5,Nframe=y,loopdiam=8}
\gasset{Nw=10,Nh=10,Nmr=5,Nframe=n}
\node(011)(35,30){\small$q$}
\node(010)(35,00){$\emptyset$}
\node(001)(10,40){\small$qr$}
\node[Nframe=y](000)(10,10){\small$r$}
\node(111)(50,45){\small$pq$}
\node[Nframe=y](110)(50,15){\small$p$}
\node[Nframe=y](101)(25,55){\small$pqr$}
\node(100)(25,25){\small$pr$}
\drawedge(000,001){}\drawedge(010,011){}
\drawedge(000,010){}\drawedge(001,011){}
\drawedge(100,101){}\drawedge(110,111){}
\drawedge(100,110){}\drawedge(101,111){}
\drawedge(100,000){}\drawedge(101,001){}
\drawedge(110,010){}\drawedge(111,011){}
\end{picture}}
\node[Nframe=n,Nw=50,Nh=60](ZES)(90,30){\unitlength 0.50mm
\begin{picture}(50,60)(5,0)
\gasset{AHnb=0}
\gasset{Nw=10,Nh=10,Nmr=5,Nframe=n}
\node(011)(35,30){\small$q$}
\node(010)(35,00){$\emptyset$}
\node[Nframe=y](001)(10,40){\small$qr$}
\node[Nframe=y](000)(10,10){\small$r$}
\node[Nframe=y](111)(50,45){\small$pq$}
\node[Nframe=y](110)(50,15){\small$p$}
\node[Nframe=y](101)(25,55){\small$pqr$}
\node(100)(25,25){\small$pr$}
\drawedge(000,001){}\drawedge(010,011){}
\drawedge(000,010){}\drawedge(001,011){}
\drawedge(100,101){}\drawedge(110,111){}
\drawedge(100,110){}\drawedge(101,111){}
\drawedge(100,000){}\drawedge(101,001){}
\drawedge(110,010){}\drawedge(111,011){}
\end{picture}}
\gasset{AHnb=1,ATnb=1,linewidth=1.3,linecolor=Gray,ATangle=30,ATLength=6,ATlength=3,AHangle=30,AHLength=6,AHlength=3}
\drawedge(EEN,TWEE){\large$+V$}
\drawedge(TWEE,DRIE){\large$*V$}
\drawedge(DRIE,VIER){\large$+V$}
\drawedge(VIER,VIJF){\large$*V$}
\drawedge(VIJF,ZES){\large$+V$}
\drawedge(ZES,EEN){\large$*V$}
\end{picture}}%
%
%
}
\caption{The orbit of a set system under vertex flips on the ground set, cf.\ Example~\ref{ex:setsystem-orbit}.}\label{fig:setsystem-orbit}
\end{figure}
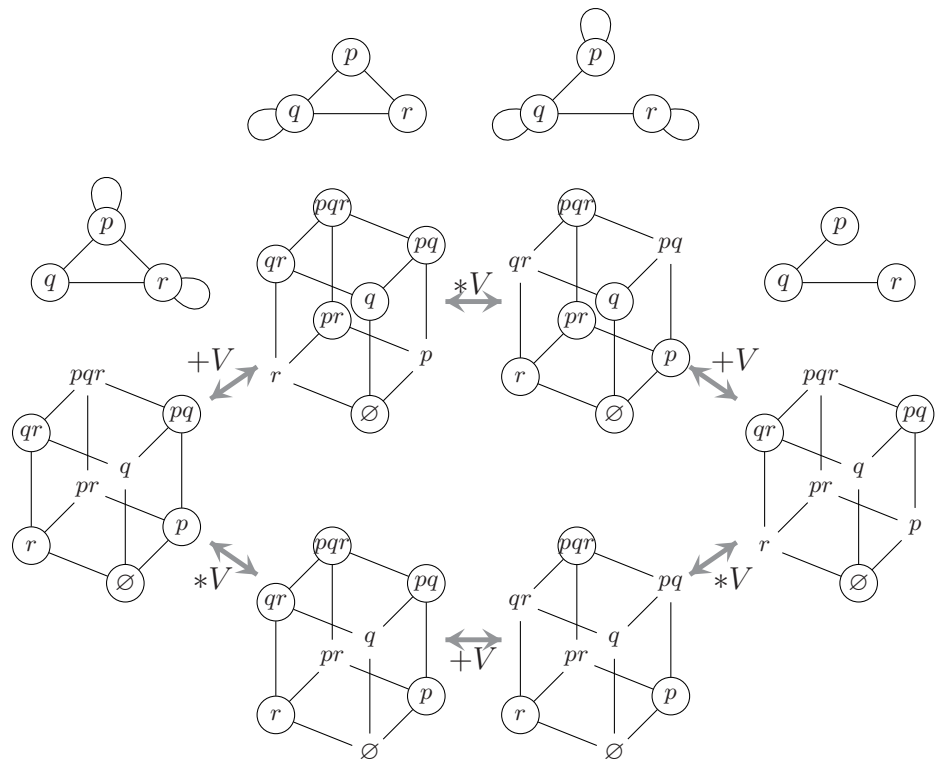

\begin{Example}\label{ex:setsystem-orbit}
Let $V = \{p,q,r\}$. Figure~\ref{fig:setsystem-orbit}
shows the set systems in the orbit under vertex
flips on $V$ of the leftmost set system $M = (\; V, \{\; \{p,q\}, \{q,r\}, \{p\},
\{r\}, \emptyset \;\} \;)$.
The set systems are depicted as a Hasse diagram where the sets that belong
to the set system are indicated by a circle (where, e.g., ``$qr$'' denotes the set $\{q,r\}$). In Section~\ref{sec:graphs} we learn that
the topmost four set systems
represent graphs, and these graphs are also indicated in the figure.
\end{Example}

It is shown in \cite{BH/PivotLoopCompl/09} that vertex flips and the removal of elements from the ground set commute when applied on different elements. Moreover,
$M+u \vertexrem u = M \vertexrem u$. This is explicitly stated as a
lemma.
\begin{Lemma}[\cite{BH/PivotLoopCompl/09}] \label{lem:comm_prop_rm_vflip}
Let $M$ be a set system and $u,v \in V$ with $u \not=v$. Then
$M+u\vertexrem v = M\vertexrem v+u$, $M*u\vertexrem v = M\vertexrem
v*u$, and $M+u \vertexrem u = M \vertexrem u$.
\end{Lemma}

Assume that $M$ is nonempty. For $X \subseteq V$, we define
$d_M(X) = \min(\{|X \sdif Y| \mid Y \in M\})$. Hence, $d_M(X)$
is the smallest distance between $X$ and the sets of $M$, where
the distance between two sets is defined as the number of
elements in the symmetric difference. We set $d_M =
d_M(\emptyset)$, the cardinality of a smallest set in $M$.

\begin{Lemma}[\cite{BH/NullityLoopcGraphsDM/10}] \label{lem:ss_distance_prop}
Let $M$ be a nonempty set system. Then $d_{M*Z}(X) = d_{M}(X
\sdif Z)$ for all $X,Z \subseteq V$, and $d_{M+Y} = d_M$ for
all $Y \subseteq V$.
\end{Lemma}

In particular, we have by Lemma~\ref{lem:ss_distance_prop}, $d_{M*Z}
= d_{M}(Z)$ for all $Z \subseteq V$.

A \emph{\dmatroid} is a nonempty set system $M$ that satisfies
the \emph{symmetric exchange axiom}: For all $X,Y \in M$ and
all $u \in X \sdif Y$, either $X \sdif \{u\} \in M$ or there is
a $v \in X \sdif Y$ with $v \not= u$ such that $X \sdif \{u,v\}
\in M$ \cite{bouchet1987}. Note that \dmatroids are closed
under pivot, i.e., $M*X$ for $X \subseteq V$ is a \dmatroid
when $M$ is a \dmatroid. Also, if $M \vertexrem u$ with $u \in
V$ is nonempty, then $M \vertexrem u$ is a \dmatroid when $M$ is a \dmatroid. A \dmatroid $M$ is called \emph{even} if the cardinalities of all sets in $M$ are of equal parity.

It turns out that a set system is an
equicardinal \dmatroid iff it is a matroid described by its family of bases
\cite[Proposition~3]{DBLP:conf/ipco/Bouchet95}. In this way,
the notion of \dmatroid is a generalization of the notion of
matroid. Note that if $M$ is a matroid (described by its bases), then $d_M$ and
$d_{M*V} = d_{M^*}$ are the rank and nullity of $M$. We implicitly assume that a matroid $M$ is described by its bases when viewing $M$ as a \dmatroid.

The result of applying loop complementation to a \dmatroid is not necessarily a \dmatroid \cite[Example~10]{BH/PivotLoopCompl/09}. We say that a \dmatroid $M$ is \emph{\vfsafe} if for any
sequence $\varphi$ of vertex flips (equivalently, pivots and
loop complementations) over $V$ we have that $M\varphi$ is a
\dmatroid. We say that a matroid $M$ is
\emph{\vfsafe} if $M$ is \vfsafe as a \dmatroid.
The class of \vfsafe \dmatroids
strictly contains the class of quaternary matroids,
i.e., matroids representable over $\gffour$ \cite{BH/QuaternaryVF}.
Recall that the class of quaternary matroids in turn strictly contains the class of binary matroids.
It is conjectured
that the class of \vfsafe matroids is equal to the
class $\mathcal{N}$ from \cite{DBLP:journals/jct/GeelenGK00},
consisting of matroids that have no minors isomorphic
to $U_{2,6}$, $U_{4,6}$, $P_6$, $F_7^-$, or $(F_7^-)^*$.
For a description of these matroids, see
\cite{Oxley/MatroidBook-2nd}.

\section{Equivalence of vf-safe delta-matroids and tight $3$-matroids}
\label{sec:equivalence_vf-safe}

In this section we show that vf-safe \dmatroids may be viewed as ``fragments'' of tight $3$-matroids. Moreover, we show that given a vf-safe \dmatroid, it is possible to reconstruct the entire tight $3$-matroid (up to renaming elements from the ground set). In this way, vf-safe \dmatroids and tight $3$-matroids turn out to be essentially equivalent.

\subsection{Unique tight extension} \label{ssec:extend_tight}
First we consider a general result that says that there is at most one way to extend a $k$-matroid to a tight $(k+1)$-matroid.

\begin{Theorem} \label{thm:extend_tight_mm}
Let $(U,\Omega)$ be a carrier. Also, let $T \in \mathcal{T}(\Omega)$ and $\Omega' = \{ \omega \setminus T \mid \omega \in \Omega\}$. If $Z'$ is a nondegenerate multimatroid with carrier $(U\setminus T,\Omega')$, then there is at most one tight multimatroid $Z$ with carrier $(U,\Omega)$ such that $Z[U\setminus T] = Z'$.
\end{Theorem}
\begin{QProof}
Let $Z_1$ and $Z_2$ be tight multimatroids such that $Z_1[U\setminus T] = Z' = Z_2[U\setminus T]$ is nondegenerate. Let $X \in \mathcal{T}(\Omega)$. We show by induction on $|X \cap T|$ that $n_{Z_1}(X) = n_{Z_2}(X)$ (this is sufficient as $Z'$ is nondegenerate). If $|X \cap T| = 0$, then $X \in \mathcal{T}(\Omega')$ and thus $n_{Z_1}(X) = n_{Z'}(X) = n_{Z_2}(X)$. Assume that the assertion holds for all $X \in \mathcal{T}(\Omega)$ with $|X \cap T| = n$. Let $|X \cap T| = n+1$. Let $x \in X \cap T$ belong to the skew class $\omega \in \Omega$. Let $p$ be any skew pair with $x \in p$. Then $|(X \sdif p) \cap T| = n$, and thus, by the induction hypothesis, $n_{Z_1}(X \sdif p) = n_{Z_2}(X \sdif p)$. Hence, if we consider $S = X \setminus \{x\}$, then $n_{Z_1}(S \cup \{y\}) = n_{Z_2}(S \cup \{y\})$ for all $y \in \omega \setminus \{x\}$. Since $Z_1$ and $Z_2$ are tight, the value of $n_{Z_1}(S \cup \{x\})$ (and $n_{Z_2}(S \cup \{x\})$) is uniquely determined given $n_{Z_1}(S \cup \{y\}) = n_{Z_2}(S \cup \{y\})$ for all $y \in \omega \setminus \{x\}$. Thus, $n_{Z_1}(X) = n_{Z_1}(S \cup \{x\}) = n_{Z_2}(S \cup \{x\}) = n_{Z_2}(X)$.
\end{QProof}

Theorem~\ref{thm:extend_tight_mm} raises the problem of characterizing those $k$-matroids $Z$ that are extendible to a tight $(k+1)$-matroid. We consider the case $k=2$ in Subsections~\ref{ssec:vfsafe_from_tight3} and \ref{ssec:tight3_from_vfsafe}.

\subsection{Tight $2$-matroids and even \dmatroids} \label{ssec:tight_2m_even_dm}

In this subsection we recall from \cite{DBLP:journals/siamdm/Bouchet97,DBLP:journals/ejc/Bouchet01} that $2$-matroids correspond to \dmatroids, and tight $2$-matroids correspond to even \dmatroids.

An $(\ell,k)$-\emph{carrier}, for positive integers $\ell$ and $k$, is a carrier $(U,\Omega)$ such that $|\Omega| = \ell$ and for all $\omega \in \Omega$, $|\omega|=k$. A \emph{transversal $k$-tuple} is a sequence $(T_1,\ldots,T_k)$ of mutually disjoint transversals of $\Omega$.
A \emph{projection} of $(U,\Omega)$ is a surjective function $\pi: U
\rightarrow V$ such that $\pi(x)=\pi(y)$ iff $x$ and $y$ are in the
same skew class $\omega \in \Omega$. Thus each skew class is assigned
by $\pi$ to a unique element of $V$.

Let $Z$ be a $2$-matroid with carrier $(U,\Omega)$,
and fix a projection $\pi: U \rightarrow V$ of $(U,\Omega)$.
Let $\tau = (T_1,T_2)$ be a transversal $2$-tuple of $\Omega$. Then the \emph{section} of $Z$ by $\tau$, denoted by $\mathcal{M}_{Z,\tau, \pi}$, is the set system $(V,D)$ with $D = \{ \pi(X \cap T_2) \mid X \in \mathcal{B}_Z \}$. If $\pi$ is clear from the context, then we often write simply $\mathcal{M}_{Z,\tau}$ to denote $\mathcal{M}_{Z,\tau, \pi}$. 
A section of a $2$-matroid is a \dmatroid, see \cite[Proposition~4.2]{DBLP:journals/siamdm/Bouchet97}. We define, for all $v \in V$, $\tau*v = (T_1 \sdif p, T_2 \sdif p)$ where $p = \pi^{-1}(v)$. It is easy to verify that $\mathcal{M}_{Z,\tau*v}*v = \mathcal{M}_{Z,\tau}$ (or see part of the proof of Theorem~\ref{thm:section_3m_vfsafe}).

Conversely, given a \dmatroid $M$ over $V$, a $(|V|,2)$-carrier $(U,\Omega)$, a transversal $2$-tuple $\tau$ of $\Omega$, and a projection $\pi: U \rightarrow V$, we may construct the $2$-matroid $Z$, denoted by $\mathcal{Z}_{M,\tau}$, such that $M = \mathcal{M}_{Z,\tau}$. We have that $\mathcal{B}_{Z} = \{ X \in \mathcal{T}(\Omega) \mid \pi(X \cap T_2) \in M \}$. The $2$-matroid $\mathcal{Z}_{M,\tau}$ is called the \emph{lift} of $M$ with respect to $\tau$ \cite[Construction~3.5]{DBLP:journals/ejc/Bouchet01}. Hence a section of a $2$-matroid $Z$ retains all essential information of $Z$.

It is shown in \cite[Theorem~5.3]{DBLP:journals/ejc/Bouchet01} that $Z$ is tight iff there is a  transversal $2$-tuple $\tau$ of $Z$ such that \dmatroid $\mathcal{M}_{Z,\tau}$ is even iff for all transversal $2$-tuples $\tau$ of $Z$, $\mathcal{M}_{Z,\tau}$ is even.

\subsection{vf-safe delta-matroids from tight $3$-matroids} \label{ssec:vfsafe_from_tight3}
In this subsection and the next we characterize the family of \dmatroids for which the corresponding $2$-matroid is extendible (in the sense of Theorem~\ref{thm:extend_tight_mm}) to a tight $3$-matroid. The approach is to extend Subsection~\ref{ssec:tight_2m_even_dm} by including both $\dual$ and $+$ and an additional transversal $T_3$.

\smallskip
\emph{In this subsection we fix a $(|V|,3)$-carrier $(U,\Omega)$, a transversal $3$-tuple $\tau$ of $(U,\Omega)$, and a projection $\pi:U \rightarrow V$ of $(U,\Omega)$.}

\smallskip
Let $v \in V$, $\omega = \pi^{-1}(v) \in \Omega$, and $\tau = (T_1,T_2,T_3)$. Let, for $i \in \{1,2,3\}$, $p_i \subseteq \omega$ be the (unique) skew pair of $\omega$ with $p_i \cap T_i = \emptyset$.
Then we define
\[
\begin{array}{rcl}
\tau*v  &=&(\, T_1\sdif p_3,\, T_2\sdif p_3,\, T_3\, ),\\
\tau+v  &=&(\, T_1,\, T_2\sdif p_1,\, T_3\sdif p_1\, ),\\
\tau\dual v&=&(\, T_1\sdif p_2,\, T_2,\, T_3\sdif p_2\, ).
\end{array}
\]
Note that $+v$, $\dual v$, and $*v$ generates a group isomorphic to $S_3$. In particular, $\tau\dual v = ((\tau +v)*v)+v$. Moreover, these operations commute on distinct elements of $V$.
Again, we assume left associativity of these operations. We extend this notation to sets, and write, for all $L \subseteq V$, e.g., $\tau+L$ to denote applying $+v$ for all $v \in L$ (in arbitrary order).

Let $Z$ be a $3$-matroid over $(U,\Omega)$.
Since $Z[T_1 \cup T_2]$ is a $2$-matroid, we have by Subsection~\ref{ssec:tight_2m_even_dm} that the section $\mathcal{M}_{Z[T_1 \cup T_2],(T_1,T_2),\pi}$ is a \dmatroid. We denote $\mathcal{M}_{Z[T_1 \cup T_2],(T_1,T_2),\pi}$ simply by $\mathcal{M}_{Z,\tau,\pi}$. Again we drop the subscript $\pi$ when $\pi$ is clear from the context. We now show that $\mathcal{M}_{Z,\tau}$ is vf-safe when $Z$ is tight.

\begin{Theorem} \label{thm:section_3m_vfsafe}
Let $Z$ be a tight $3$-matroid over $(U,\Omega)$.
For $v \in V$, $\mathcal{M}_{Z,\tau} = \mathcal{M}_{Z,\tau+v}+v = \mathcal{M}_{Z,\tau*v}*v = \mathcal{M}_{Z,\tau \dual v}\dual v$. In particular, $\mathcal{M}_{Z,\tau}$ is vf-safe.
\end{Theorem}
\begin{QProof}
Let $\tau = (T_1,T_2,T_3)$.

For all $Y\subseteq V$, there is a unique ``lift'' $X \in \mathcal{T}(\Omega)$ such that $\pi(X\cap T_2) = Y$ and $X \cap T_3 = \emptyset$, which we denote by $Y^{(\tau)}$. Thus, $Y \in \mathcal{M}_{Z,\tau}$ iff $Y^{(\tau)} \in \mathcal{B}_Z$.

Let, for $i \in \{1,2,3\}$, $p_i$ be the (unique) skew pair of $\pi^{-1}(v)$ with $p_i \cap T_i = \emptyset$.

We first consider $*$. Note that $\tau*v=(T_1\sdif p_3,T_2\sdif p_3,T_3)$. Hence $Y^{(\tau*v)} = (Y \sdif \{v\})^{(\tau)}$. We have $Y\in \mathcal{M}_{Z,\tau}*v$ iff $Y \sdif \{v\} \in \mathcal{M}_{Z,\tau}$ iff $(Y \sdif \{v\})^{(\tau)} = Y^{(\tau*v)} \in  \mathcal{B}_Z$ iff $Y\in \mathcal{M}_{Z,\tau*v}$.

We now consider $+$. Note that $\tau+v=(T_1,T_2\sdif p_1,T_3\sdif p_1)$. Assume first that $v \notin Y$. Then $Y^{(\tau)} = Y^{(\tau+v)}$ as $Y^{(\tau)} \cap \omega \subseteq T_1$. We have $Y\in \mathcal{M}_{Z,\tau}+v$ iff $Y\in \mathcal{M}_{Z,\tau}$ iff $Y^{(\tau)} = Y^{(\tau+v)} \in \mathcal{B}_Z$ iff $Y\in \mathcal{M}_{Z,\tau+v}$.

Assume now that $v\in Y$. Then $Y^{(\tau)} = (Y \setminus \{v\})^{(\tau+v)}$. We have $Y \in \mathcal{M}_{Z,\tau} + v$ iff either $Y \in \mathcal{M}_{Z,\tau}$ or $Y \setminus \{v\} \in \mathcal{M}_{Z,\tau}$ but not both iff either $Y^{(\tau)} \in \mathcal{B}_Z$ or $(Y \setminus \{v\})^{(\tau)} \in \mathcal{B}_Z$ but not both iff $Y^{(\tau)} \sdif p_1 \in \mathcal{B}_Z$ (since $Z$ is tight) iff $Y^{(\tau+v)} \in \mathcal{B}_Z$ iff $Y \in \mathcal{M}_{Z,\tau+v}$.

We finally consider $\dual$, which is now easy. We have $\mathcal{M}_{Z,\tau}\dual v = \mathcal{M}_{Z,\tau}+v*v+v = \mathcal{M}_{Z,\tau+v*v+v} = \mathcal{M}_{Z,\tau\dual v}$.
\end{QProof}

\subsection{Tight $3$-matroids from vf-safe delta-matroids}
\label{ssec:tight3_from_vfsafe}

A \emph{pre-multimatroid} $Z$ is a triple $(U,\Omega,\mathcal{I})$, where $(U,\Omega)$ is a carrier and $\mathcal{I} \subseteq \mathcal{S}(\Omega)$. Note that we may define, e.g., the restriction $Z[X]$ for $X \in \mathcal{S}(\Omega)$ exactly as we did for multimatroids.

\smallskip
\emph{In this subsection we fix a nonempty set system $M$ over $V$, a $(|V|,3)$-carrier $(U,\Omega)$, a transversal $3$-tuple $\tau = (T_1,T_2,T_3)$ of $\Omega$, and a projection $\pi:U \rightarrow V$ of $(U,\Omega)$.}

\smallskip

We define
$$
\mathcal{B}_{M,\tau, \pi} = \{\; X \in \mathcal{T}(\Omega) \mid
     \pi(X \cap T_2) \in M \dual \pi(X \cap T_3) \;\}.
$$
For symmetry's sake,
note that the condition in the definition of $\mathcal{B}_{M,\tau,\pi}$ is equivalent to $\emptyset \in M \dual \pi(X \cap T_3) * \pi(X \cap T_2) + \pi(X \cap T_1)$ by Lemma~\ref{lem:ss_distance_prop}.
As the sets $\pi(X\cap T_i)$, $i \in \{1,2,3\}$, are disjoint, we can perform these operations in any order.

We denote by $\mathcal{Z}_{M,\tau,\pi}$ the pre-multimatroid $(U,\Omega,\mathcal{B}_{M,\tau,\pi})$. Again, if $\pi$ is clear from the context, then we often write simply $\mathcal{B}_{M,\tau}$ and $\mathcal{Z}_{M,\tau}$ to denote $\mathcal{B}_{M,\tau,\pi}$ and $\mathcal{Z}_{M,\tau,\pi}$, respectively.

Note that if $M$ is a \dmatroid, then, by Subsection~\ref{ssec:tight_2m_even_dm}, $\mathcal{Z}_{M,\tau}[T_1 \cup T_2]$ is a $2$-matroid and the lift of $M$ with respect to $(T_1,T_2)$. Thus, $\mathcal{Z}_{M,\tau}[T_1 \cup T_2] = \mathcal{Z}_{M,(T_1,T_2)}$.

\begin{Lemma}
\label{lem:single_el_exchg}
Let $v \in V$. Then $\mathcal{Z}_{M,\tau} = \mathcal{Z}_{M+v,\tau+v} = \mathcal{Z}_{M\dual v,\tau \dual v} = \mathcal{Z}_{M* v,\tau * v}$
\end{Lemma}
\begin{QProof}
For later use, we first show a general property of vertex flip operations on a set system $M$.
Let $B,C \subseteq V$ be disjoint, and let $Y\subseteq V$.
Then
\begin{eqnarray}
M+Y*B\dual C &=& M  *(B\sdif Y') \dual {}(C\sdif Y') +Y
 \mbox{ with } Y'= Y \cap (B\cup C), \label{eqn:lemma_loopc}\\
M * Y * B \dual C &=& M * (B\sdif Y') \dual C  + (Y\cap C)
 \mbox{ with } Y'= Y \setminus C, \label{eqn:lemma_pivot}\\
M \dual Y * B \dual C &=& M * B \dual {}(C \sdif Y') + (Y\cap B)
 \mbox{ with } Y'= Y \setminus B. \label{eqn:lemma_dpivot}
\end{eqnarray}

As the vertex flips on different vertices commute, it suffices to consider only a single element $v \in V$. In other words, we may assume without loss of generality that $Y$, $B$, and $C$ are all subsets of $\{v\}$.
If $v \notin Y$, then there is nothing to prove. Assume $v \in Y$.
Then, depending on whether $v$ is in $B$ or $C$, we observe that
$ M + v * v = M \dual v + v $
(for $v\in B$),
$ M + v \dual v = M * v + v $
(for $v\in C$), and
$ M + v  = M + v $
(for $v\notin B\cup C$). This proves Equality~(\ref{eqn:lemma_loopc}).

Similarly,
$ M * v * v = M  $
(for $v\in B$),
$ M * v \dual v = M \dual v + v $
(for $v\in C$), and
$ M * v  = M  * v $
(for $v\notin B\cup C$), which leads to Equality~(\ref{eqn:lemma_pivot}).
Equality~(\ref{eqn:lemma_dpivot}) is proved in an analogous way.

By Lemma~\ref{lem:ss_distance_prop},
\begin{eqnarray*}
\emptyset \in M+Y*B\dual C &\mbox{ iff }& \emptyset \in M  *(B\sdif Y') \dual {}(C\sdif Y') \mbox{ with } Y'= Y \cap (B\cup C),\\
\emptyset \in M * Y * B \dual C &\mbox{ iff }& \emptyset \in M * (B\sdif Y') \dual C \mbox{ with } Y'= Y \setminus C, \\
\emptyset \in M \dual Y * B \dual C &\mbox{ iff }& \emptyset \in M * B \dual {}(C \sdif Y') \mbox{ with } Y'= Y \setminus B.
\end{eqnarray*}

Let, for $i \in \{1,2,3\}$, $p_i \subseteq \omega$ be the (unique) skew pair of $\omega$ with $p_i \cap T_i = \emptyset$. Let $X \in \mathcal{T}(\Omega)$. We thus have $X \in \mathcal{B}_{M+v,\tau}$ iff $\emptyset \in M+v * \pi(X \cap T_2) \dual \pi(X \cap T_3)$ iff $\emptyset \in M * (\pi(X \cap T_2) \sdif Y') \dual (\pi(X \cap T_3) \sdif Y')$ with $Y' = \{v\} \cap (\pi(X \cap T_2) \cup \pi(X \cap T_3)) = \{v\} \setminus \pi(X \cap T_1)$. We have that $\pi(X \cap (T_2 \sdif p_1))$ is equal to $\pi(X \cap T_2)$ if $v \in \pi(X \cap T_1)$ and equal to $\pi(X \cap T_2) \sdif \{v\}$ otherwise, and similarly for $\pi(X \cap (T_3 \sdif p_1))$. Hence $\emptyset \in M * (\pi(X \cap T_2) \sdif Y') \dual (\pi(X \cap T_3) \sdif Y')$ iff $\emptyset \in M * \pi(X \cap (T_2 \sdif p_1)) \dual \pi(X \cap (T_3 \sdif p_1))$ iff $X \in \mathcal{B}_{M,\tau+v}$.

Similarly, we obtain $X \in \mathcal{B}_{M*v,\tau}$ iff $X \in \mathcal{B}_{M,\tau*v}$ and $X \in \mathcal{B}_{M\dual v,\tau}$ iff $X \in \mathcal{B}_{M,\tau+v}$.
\end{QProof}

We now show that $\mathcal{Z}_{M,\tau}$ is a tight 3-matroid when $M$ is a vf-safe \dmatroid.
\begin{Theorem} \label{thm:vfsafe_dmatroid_tight3matroid}
Let $M$ be a vf-safe \dmatroid. Then $\mathcal{Z}_{M,\tau}$ is a tight $3$-matroid described by its bases.
\end{Theorem}
\begin{QProof}
Let $T \in \mathcal{T}(\Omega)$.
We first argue that $Z' = \mathcal{Z}_{M,\tau}[T]$ is a matroid.
Let $\varphi$ be a sequence of $+$, $*$, and $\dual$ operations such that $\tau\varphi = (T'_1,T'_2,T'_3)$ with $T \subseteq T'_1\cup T'_2$.
Then by Lemma~\ref{lem:single_el_exchg}, $\mathcal{Z}_{M,\tau} = \mathcal{Z}_{M\varphi,\tau\varphi}$. By Subsection~\ref{ssec:tight_2m_even_dm},
$\mathcal{Z}_{M\varphi,\tau\varphi}[T_1\cup T_2] = \mathcal{Z}_{M\varphi,(T'_1,T'_2)}$ is the lift of \dmatroid $M\varphi$ with respect to $(T'_1,T'_2)$. Consequently, $Z' = \mathcal{Z}_{M\varphi,(T'_1,T'_2)}[T]$ is a matroid.

We now show the second defining property of multimatroids (in Definition~\ref{def:multimatroid}). Let $\mathcal{I} = \{ I \subseteq X \mid X \in \mathcal{B}_{M,\tau} \}$. Let $I \in \mathcal{I}$ and $\omega \in \Omega$ with $\omega\cap I = \emptyset$, and let $X \in \mathcal{B}_{M,\tau}$ with $I \subseteq X$. Consider now $Y=X\setminus \omega$. Let $M'= M * \pi(Y \cap T_2) \dual \pi(Y \cap T_3)$. We have (\emph{i}) $Y\cup (\omega \cap T_1) \in \mathcal{B}_{M,\tau}$ iff $\emptyset \in M'$, (\emph{ii}) $Y\cup (\omega \cap T_2) \in \mathcal{B}_{M,\tau}$ iff $\emptyset \in M'*v$ iff $\{v\} \in M'$, and (\emph{iii}) $Y\cup (\omega \cap T_3) \in \mathcal{B}_{M,\tau}$ iff $\emptyset \in M' \dual v$ iff exactly one of $\emptyset$ and $\{v\}$ belongs to $M'$. Thus exactly two of the transversals $Y\cup \{u\}$, $u \in \omega$, belong to $\mathcal{B}_{M,\tau}$. Consequently, any skew pair $p\subseteq \omega$ intersects with at least one of these two transversals in $\mathcal{B}_{M,\tau}$. Hence $\mathcal{Z}_{M,\tau}$ is a $3$-matroid described by its bases. Moreover, by Proposition~\ref{prop:base-tight}, the 3-matroid $\mathcal{Z}_{M,\tau}$ is tight.
\end{QProof}

Recall from Subsection~\ref{ssec:tight_2m_even_dm} that the lift $M \mapsto\mathcal{Z}_{M,(T_1,T_2)}$ is a one-to-one correspondence from the family of \dmatroids over $V$ to the family of $2$-matroids with carrier $(T_1\cup T_2,\{\omega\setminus T_3 \mid \omega \in \Omega\})$, with inverse mapping $Z\mapsto \mathcal{M}_{Z,(T_1,T_2)}$.

By Theorems~\ref{thm:section_3m_vfsafe} and \ref{thm:vfsafe_dmatroid_tight3matroid}, for a vf-safe \dmatroid $M$ over $V$, the $3$-matroid $\mathcal{Z}_{M,\tau}$ is the (unique) extension of the $2$-matroid $\mathcal{M}_{Z,(T_1,T_2)}$ to a tight $3$-matroid over $(U,\Omega)$ (as in Theorem~\ref{thm:extend_tight_mm}). The $2$-matroids involved are precisely the $2$-matroids that allow for such a (unique) extension.

Let $\Theta$ be the mapping from $3$-matroids over $(U,\Omega)$ to $2$-matroids which sends each $Z$ to the restriction $Z[T_1\cup T_2]$. Then, for tight $3$-matroids over $(U,\Omega)$, $Z \mapsto\mathcal{M}_{Z,\tau} = \mathcal{M}_{\Theta(Z),(T_1,T_2)}$ is the inverse of $M \mapsto\mathcal{Z}_{M,\tau}$.

This leads to the following result.
\begin{Theorem} \label{thm:vfsafe_equiv_t3matroid}
The mapping $M \mapsto \mathcal{Z}_{M,\tau}$ is a one-to-one correspondence from the family of vf-safe \dmatroids over $V$ to the family of tight $3$-matroids over $(U,\Omega)$, with inverse $Z \mapsto\mathcal{M}_{Z,\tau}$.
\end{Theorem}

\subsection{Carrying over nullity and minor}
In this subsection we tie the notion of nullity in a tight 3-matroid to that of distance in the corresponding vf-safe \dmatroid.
We also carry over the concept of minor for tight $3$-matroids to vf-safe \dmatroids
(and similarly for $2$-matroids to \dmatroids).

\begin{Lemma} \label{lem:null=dist}
Let $M$ be a vf-safe \dmatroid over $V$. Let $(U,\Omega)$ be a $(|V|,3)$-carrier, $\tau$ be a transversal $3$-tuple of $\Omega$, and $\pi:U \rightarrow V$ be a projection of $(U,\Omega)$. Then for all $T \in \mathcal{T}(\Omega)$, $n_{\mathcal{Z}_{M,\tau}}(T) = d_{M* \pi(T\cap T_2)\dual\pi(T\cap T_3)}$.
\end{Lemma}
\begin{QProof}
Let $\tau = (T_1,T_2,T_3)$ and $T \in \mathcal{T}(\Omega)$. Let $M' = M* \pi(T\cap T_2)\dual\pi(T\cap T_3)$. By 
Theorem~\ref{thm:vfsafe_dmatroid_tight3matroid},
$\mathcal{Z}_{M',\tau}$ is a $3$-matroid. By Lemma~\ref{lem:single_el_exchg}, $n_{\mathcal{Z}_{M',\tau}}(T) = n_{\mathcal{Z}_{M,\tau'}}(T)$ with $\tau' = \tau  * \pi(T\cap T_2)\dual\pi(T\cap T_3)$. Now, $\tau' = (T'_1,T'_2,T'_3)$ with $T'_1 = T$. Let $Z = \mathcal{Z}_{M,\tau'}$.

We have that $r_Z(T) = r_{Z[T\cup T'_2]}(T) = r((Z[T\cup T'_2])[T])$ which is equal to the largest independent set $I$ of $Z[T\cup T'_2]$ with $I \subseteq T$. Hence, $r_Z(T)$ is the maximum value of $|X\cap T|$ with $X \in \mathcal{B}_{Z[T\cup T'_2]}$. Consequently, $n_Z(T)$ is the minimum value of $|\Omega| - |X\cap T| = |V| - |\pi(X\cap T)|$ with $X \in \mathcal{B}_{Z[T\cup T'_2]}$.

Note that $\mathcal{Z}_{M,\tau'}[T \cup T'_2]$ is the lift of $M$ with respect to $(T,T'_2)$. Hence, for all $Y \subseteq V$, $Y \in M$ iff $Y = \pi(X\cap T'_2)$ with $X \in \mathcal{B}_{Z[T\cup T'_2]}$. Thus, $d_{M}$ is the minimum value of the $|\pi(X\cap T'_2)| = |V| - |\pi(X\cap T)|$ with $X \in \mathcal{B}_{Z[T\cup T'_2]}$. Hence $n_{\mathcal{Z}_{M',\tau}}(T) = n_Z(T) = d_{M}$. By change of variables $M:=M'$ we obtain the required result.
\end{QProof}
Let $M$ be a vf-safe \dmatroid. Since $\mathcal{Z}_{M,\tau}$ is tight, we have by Lemma~\ref{lem:null=dist} that the values $d_M$, $d_{M*v}$, and $d_{M\dual v}$ are such that precisely two of the three are equal, to say $m$, and the third is equal to $m+1$. This has also been shown in a direct way (without the use of multimatroids) in \cite{BH/NullityLoopcGraphsDM/10}.

Obviously, we may formulate the ``2-matroid version'' of Lemma~\ref{lem:null=dist}: if $Z$ is a $2$-matroid with transversal $2$-tuple $\tau = (T_1,T_2)$, then $n_{Z}(T) = d_{\mathcal{M}_{Z,\tau}* \pi(T\cap T_2)}$ for all $T \in \mathcal{T}(\Omega)$.
Note that the 2-matroid version of Lemma~\ref{lem:null=dist} is \emph{not} strictly a special case of Lemma~\ref{lem:null=dist} as not all $2$-matroids allow for an extension to tight $3$-matroids (or, equivalently, not every \dmatroid is vf-safe).

We now turn to minors.
\begin{Lemma} \label{lem:minors_dm_mm}
Let $M$ be a vf-safe \dmatroid over $V$. Let $Z = \mathcal{Z}_{M,\tau}$ with $\tau = (T_1,T_2,T_3)$ a transversal $3$-tuple of $\Omega$, $(U,\Omega)$ a $(|V|,3)$-carrier, and $\pi:U \rightarrow V$ a projection of $(U,\Omega)$.

Moreover, let $v \in V$, $M_1 = M \vertexrem v$, $M_2  = M*v \vertexrem v$, $M_3 = M\dual v \vertexrem v$, and $\omega = \pi^{-1}(v)$.

For all $i \in \{1,2,3\}$, $M_i$ is nonempty iff the (unique) element $t_i \in T_i \cap \omega$ is nonsingular in $Z$. Moreover, if $t_i$ is nonsingular in $Z$, then $Z|t_i = \mathcal{Z}_{M_i,\tau'}$ with $\tau' = (T_1 \setminus \omega, T_2 \setminus \omega, T_3 \setminus \omega)$ and so $M_i$ is a vf-safe \dmatroid.
\end{Lemma}
\begin{QProof}
First, $t_1$ is nonsingular in $Z$ iff $Z[T_1\cup T_2] = \mathcal{Z}_{M,(T_1,T_2)}$ has a basis containing $t_1$ iff $M$ contains a set without $v$ iff $M\vertexrem v$ is nonempty. By Lemma~\ref{lem:single_el_exchg}, $t_2$ is nonsingular in $Z$ iff $t_1$ is nonsingular in $\mathcal{Z}_{M,\tau*v} = \mathcal{Z}_{M*v,\tau}$ iff $M*v\vertexrem v$ is nonempty. In the same way, we obtain that $t_3$ is nonsingular in $Z$ iff $M\dual v\vertexrem v$ is nonempty.

Let $t_i$ be nonsingular in $Z$ for some $i \in \{1,2,3\}$, and let $X \in \mathcal{T}(\Omega\setminus \{\omega\})$. Then, by definition of minor, $X$ is a basis of $Z|t_i$ iff $X \cup \{t_i\}$ is a basis of $Z$. Also, $X$ is a basis of the pre-multimatroid $\mathcal{Z}_{M_i,\tau'}$ iff $\emptyset \in M_i * \pi(X \cap (T_2 \setminus \omega)) \dual \pi(X \cap (T_3 \setminus \omega))$ iff $\emptyset \in M * \pi((X \cup \{t_i\}) \cap T_2) \dual \pi((X \cup \{t_i\}) \cap T_3) \vertexrem v$ iff $X \cup \{t_i\}$ is a basis of $\mathcal{Z}_{M,\tau} = Z$. Thus $Z|t_i = \mathcal{Z}_{M_i,\tau'}$. Since $Z|t_i$ is a tight $3$-matroid, $M_i$ is a vf-safe \dmatroid by Theorem~\ref{thm:vfsafe_equiv_t3matroid}.
\end{QProof}

Again we may formulate the obvious ``$2$-matroid version'' of Lemma~\ref{lem:minors_dm_mm}. This leads to the following concepts. We call $v \in V$ \emph{nonsingular} in a nonempty set system $M$ if both $M\vertexrem v$ and $M*v\vertexrem v$ are nonempty. Thus, by the ``2-matroid version'' of Lemma~\ref{lem:minors_dm_mm}, a $v$ is nonsingular in a \dmatroid $M$ iff all $u \in \pi^{-1}(v)$ are nonsingular in $\mathcal{Z}_{M,\tau}$ for suitable transversal $2$-tuple $\tau$ iff $\pi^{-1}(v)$ is nonsingular in $\mathcal{Z}_{M,\tau}$.

Note that $v \in V$ is nonsingular in a nonempty set system $M$ iff there are $X_1,X_2 \in M$ with $v \in X_1 \sdif X_2$. Clearly, for all $w \in V$, $v$ is nonsingular in $M$ iff $v$ is nonsingular in $M*w$. Note that every $v \in V$ is singular in $M$ iff $M$ contains only one subset.

Note that if $M$ is a matroid, then $M\vertexrem v$ is nonempty iff $v$ is not a coloop of $M$, and $M*v\vertexrem v$ is nonempty iff $v$ is not a loop of $M$. Hence $v$ is nonsingular in $M$ iff $v$ is neither a coloop nor a loop.

We call $v \in V$ \emph{strongly nonsingular} in a nonempty set system $M$ if
$M\vertexrem v$, $M*v\vertexrem v$, and $M\dual v\vertexrem v$
are all nonempty. Thus, by Lemma~\ref{lem:minors_dm_mm}, $v$ is strongly nonsingular in a vf-safe \dmatroid $M$ iff $\pi^{-1}(v)$ is nonsingular in $\mathcal{Z}_{M,\tau}$ for suitable transversal $3$-tuple $\tau$.

\begin{Remark}
In line with matroid theory, one may define an element $v$ of a \dmatroid $M$ to be a \emph{coloop} of $M$ if $M\vertexrem v$ is empty, and a \emph{loop} of $M$ if $M*v\vertexrem v$ is empty. Then one may define the \emph{deletion} of $v$ in $M$, here denoted by $M \bbslash v$, to be $M \vertexrem v$ if $v$ is not a coloop of $M$ and $M*v\vertexrem v$ otherwise. Similarly, one may define the \emph{contraction} of $v$ in $M$, denoted by $M \slash v$, to be $M*v\vertexrem v$ if $v$ is not a loop of $M$ and $M\vertexrem v$ otherwise. Note that if $M$ is a matroid described by its bases, then these notions of coloop, loop, deletion and contraction coincide with the usual matroid-theoretical definitions (note, e.g, that if $v$ is not a loop, then $M \slash v = (V \setminus \{v\}, \{X \setminus \{v\} \mid v \in X \in M\}) = M*v\vertexrem v$). Now, by the $2$-matroid version of Lemma~\ref{lem:minors_dm_mm}, Proposition~\ref{prop:minor_singular}, and the fact that minors on distinct elements in a multimatroid commute, we have that deletion and contraction on \dmatroids $M$ commute on distinct elements. Note that these operations do not commute for set systems $M$ in general: for example, if $M = (\{u,v,w\},\{\{u\}, \{v\}, \{v,w\}\})$, then $M\bbslash u \bbslash v = (\{w\},\{\emptyset, \{w\}\})$ and $M\bbslash v \bbslash u = (\{w\},\{\emptyset\})$. A similar remark can be made for vf-safe \dmatroids and its three types of minors using Lemma~\ref{lem:minors_dm_mm}. For notational convenience, we stick with the set system operation $\vertexrem$ in this paper, but the reader may trivially reformulate the results in this paper in terms of the minor definitions given here.
\end{Remark}

\section{Transition Polynomial for Set Systems}
\label{sec:def_sets_interlacep}

\subsection{The transition polynomial}

We move from multimatroids to set systems, and consider a, rather generic,
weighted polynomial for set systems. We first obtain a technical
result that shows how the variables of the polynomial change when
one of the vertex flip operations ${}+Y$, ${}\dual Y$ and ${}*Y$ is
applied. In the next subsections more interesting polynomials appear as
specializations.

Let $V$ be a finite set. We define $\mathcal{P}_3(V)$ to be the set
of triples $(V_1,V_2,V_3)$ where $V_1$, $V_2$, and $V_3$ are
pairwise disjoint subsets of $V$ such that $V_1 \cup V_2 \cup V_3 =
V$. Therefore $V_1$, $V_2$, and $V_3$ form an ``ordered partition''
of $V$ where $V_i = \emptyset$ for some $i \in \{1,2,3\}$ is
allowed.

We now define a weighted polynomial for set systems.
\begin{Definition}
Let $M$ be a nonempty set system over $V$. We define the \emph{(weighted)
transition polynomial} of $M$ as follows:
\begin{eqnarray*}
Q(M)(\vec{a},\vec{b},\vec{c},y) = \sum_{(A,B,C) \in \mathcal{P}_3(V)} a_A b_B c_C
\, y^{d_{M*B\dual C}},
\end{eqnarray*}
where $\vec{a}$ is a vector indexed by $V$ with entries $a_{v}$ for all $v \in V$, $a_A = \prod_{v \in A} a_{v}$, and similarly for $\vec{b}$, $b_B$, $\vec{c}$, and $c_C$.
\end{Definition}
For notational convenience, we often omit the vectors $\vec{a},\vec{b},\vec{c}$ and write simply $Q(M)(y)$. In fact, we often also omit the variable $y$ in $Q(M)(y)$ and write simply $Q(M)$.

Note that one may write $y^{d_{M+A*B\dual C}}$ instead of
$y^{d_{M*B\dual C}}$ in the definition of $Q(M)$ (indeed, as
$A$ is disjoint from $B$ and $C$, ${}+A$ commutes with the
other operations and we apply Lemma~\ref{lem:ss_distance_prop}).

In Section~\ref{sec:graph_poly} we will find that $Q(M)$
generalizes the interlace polynomial for graphs \cite{ArratiaBS04}.

\medskip
The next result shows that the vertex flip operations result in a
permutation of the variables $a_A$, $b_B$, and $c_C$.

\begin{Theorem} \label{thm:multivar_int_pol}
Let $M$ be a nonempty set system and $Y \subseteq V$. Then
\begin{eqnarray*}
Q(M+Y) &=& \sum_{(A,B,C) \in \mathcal{P}_3(V)} a_A\,  b_{(B \sdif Y')}\,
c_{(C \sdif Y')}\,  y^{d_{M*B\dual C}} \mbox{ with } Y' = Y \setminus
A,
\\
Q(M\dual Y) &=& \sum_{(A,B,C) \in \mathcal{P}_3(V)} a_{(A \sdif Y')}\,
b_{B}\,  c_{(C \sdif Y')}\,  y^{d_{M*B\dual C}} \mbox{ with } Y' = Y
\setminus B, \mbox{ and}
\\
Q(M*Y) &=& \sum_{(A,B,C) \in \mathcal{P}_3(V)} a_{(A \sdif Y')}\, b_{(B
\sdif Y')}\, c_{C}\,  y^{d_{M*B\dual C}} \mbox{ with } Y' = Y \setminus
C.
\end{eqnarray*}
\end{Theorem}
\begin{QProof}
By Equality~(\ref{eqn:lemma_loopc}) in the proof of Lemma~\ref{lem:single_el_exchg} we have $d_{M+Y*B\dual C} = d_{M  *(B\sdif Y') \dual {}(C\sdif Y')}$. The equality for $Q(M+Y)$ is obtained by changing variables $B:= B\sdif Y'$ and $C:=C\sdif Y'$. The equalities for $Q(M\dual Y)$ and $Q(M*Y)$ are obtained similarly.
\end{QProof}

\subsection{Specializations of the transition polynomial} \label{ssec:variants_ss_polynomials}
In this subsection we consider four interesting special cases of
the transition polynomial $Q(M)$, which generalizes
the two types of interlace polynomials
\cite{Aigner200411,Arratia2004199} and the bracket polynomial
\cite{Traldi/Bracket1/09} for graphs. These special cases each
fulfill a particular invariance result with respect to pivot,
loop complementation, or dual pivot. These invariance results
do not hold for $Q(M)$ in general.

In this subsection we again let $M$ be a nonempty set system.
We let $Q_{[a,b,c]}(M)$ be $Q(M)$ where $a_u = a$, $b_u = b$, and $c_u
= c$ for all $u \in V$.
Thus
$$
Q_{[a,b,c]}(M)(y) = \sum_{(A,B,C) \in \mathcal{P}_3(V)} a^{|A|} b^{|B|} c^{|C|}
y^{d_{M*B\dual C}}.
$$
Note that by
Theorem~\ref{thm:multivar_int_pol}, $Q_{[a,b,c]}(M) =
Q_{[a,c,b]}(M+V) = Q_{[c,b,a]}(M\dual V) = Q_{[b,a,c]}(M*V)$.

We consider specific values for $a$, $b$, and $c$ as specializations of $Q_{[a,b,c]}(M)$.

\paragraph{The polynomial $Q_1(M)$}
Let $Q_1(M) = Q_{[1,1,1]}(M)$. Explicitly, we have
$$
Q_1(M)
 = \sum_{X,Y \subseteq V, X \cap Y = \emptyset} y^{d_{M*Y\dual X}}
 = \sum_{X \subseteq V} \sum_{Z \subseteq X} y^{d_{M+Z}(X)}.
$$
The last equality can be seen as follows.
For $X \cap Y = \emptyset$ we have
$M*Y\dual X = M*Y+X*X+X = M+X*Y*X+X = M+X*(X \cup Y)+X$.
Moreover, $d_{M+X*(X \cup Y)+X} = d_{M+X*(X \cup Y)} =
d_{M+X}(X \cup Y)$. Finally, change variables [$Z:=X$, $X:=X \cup Y$].

By Theorem~\ref{thm:multivar_int_pol}, $Q_1(M)$ is invariant under
pivot, loop complementation, and dual pivot.
We will see in Section~\ref{sec:graph_poly} that
$Q_1(M)$ generalizes an interlace polynomial for simple graphs
defined in \cite{Aigner200411}.

\paragraph{The polynomials $q_i(M)$}
We now give polynomials $q_1(M)$, $q_2(M)$, and $q_3(M)$ which are
invariant under pivot, loop complementation, and dual pivot,
respectively.
In Section~\ref{sec:graph_poly} we find that $q_1(M)$
generalizes the (single-variable) interlace polynomial (from
\cite{Arratia2004199}), and that $q_2(M)$ generalizes the bracket
polynomial for graphs defined in \cite{Traldi/Bracket1/09}.

\medskip
Let $q_1(M) = Q_{[1,1,0]}(M)$.
Thus, by Lemma~\ref{lem:ss_distance_prop},
$$
q_1(M)
= \sum_{X \subseteq V} y^{d_{M*X}}
= \sum_{X \subseteq V} y^{d_{M}(X)}.
$$
We remark that $q_1(M)$ is a very natural polynomial as it counts the distances of each $X \subseteq V$ with $M$. By Theorem~\ref{thm:multivar_int_pol} (or by the explicit
formulation above), $q_1(M)$ is invariant under pivot.

\medskip
Let $q_2(M) = q_1(M\dual V) = Q_{[0,1,1]}(M)$. Recall from above (regarding $Q_1$) that for $X,Y \subseteq V$ with $X \cap Y = \emptyset$, we have
$d_{M*Y\dual X} = d_{M+X}(X \cup Y)$. Therefore, considering the case $Y = V\setminus X$, we have
$$
q_2(M)
= \sum_{X \subseteq V} y^{d_{M*(V\setminus X)\dual X}}
= \sum_{X \subseteq V} y^{d_{M+X}(V)}.
$$
By Theorem~\ref{thm:multivar_int_pol} (or by the explicit
formulation above), $q_2(M)$ is invariant under loop
complementation.

\medskip
Finally $q_3(M) = q_1(M+V) = Q_{[1,0,1]}(M)$ completes the family.
Thus,
$$
q_3(M)
= \sum_{X \subseteq V} y^{d_{M\dual X}}
= \sum_{X \subseteq V} y^{d_{M+X}(X)}.
$$
The latter equality holds, since
$d_{M\dual X} = d_{M+X*X+X} = d_{M+X*X} = d_{M+X}(X)$.

By Theorem~\ref{thm:multivar_int_pol} (or by the explicit
formulation above), $q_3(M)$ is invariant under dual
pivot. Moreover, by Theorem~\ref{thm:multivar_int_pol}, we find
that $q_3(M) = q_2(M*V)$.
Note that since pivot, loop complementation, and dual pivot are
involutions, we also have, e.g., $q_2(M) = q_3(M*V)$.

\begin{Example}
Let $M = (\; V, \{\; \{p,q\}, \{q,r\}, \{p\}, \{r\}, \emptyset
\;\} \;)$, cf.\ Example~\ref{ex:setsystem-orbit}. Then $Q_1(M)
= 16+10y+y^2$, $q_1(M) = 5+3y$, $q_2(M) = q_1(M\dual V) =
3+4y+y^2$, and $q_3(M) = q_1(M+V) = 6+2y$.
\end{Example}

Note that the invariance properties of $Q_1(M)$ and $q_i(M)$ ($i \in
\{1,2,3\}$) do \emph{not} hold in general for the transition polynomial $Q(M)$
itself. In the next sections we will focus on particular
recursive relations for $Q_1(M)$ and $q_i(M)$ ($i \in \{1,2,3\}$), which
also do not hold for $Q(M)$ in general.

\subsection{Two-variable polynomials and the Tutte polynomial}
\label{ssec:two-var_&_tutte} We often consider single-variable polynomials $Q_1(M)$ and $q_i(M)$
by letting $a$, $b$, and $c$ be constants in $Q_{[a,b,c]}(M)$.
However, at some points it is useful to consider some of
$a$, $b$, and $c$ as variables.

Again, let $M$ be a nonempty set system. We have $Q_{[a,b,0]}(M)(y) = \sum_{X \subseteq V}
a^{|V\setminus X|}  b^{|X|} y^{d_{M*X}}$. Let $\bar q(M)(x,y) =
Q_{[1,x,0]}(M)(y)$. Then obviously $\bar q(M)(x,y) = \sum_{X
\subseteq V} x^{|X|} y^{d_{M*X}}$ and $\bar q(M)(1,y) =
q_1(M)(y)$. Note however that, unlike the single-variable case
$q_1(M)$, $\bar q(M)(x,y)$ is not invariant under pivot.
Polynomial $\bar q(M)(x,y)$ will be of interest in
Section~\ref{sec:graph_poly}. Another two-variable case of interest in
Section~\ref{sec:graph_poly} is $\bar q(M *V \dual V) =
Q_{[0,1,x]}(M)$.

\medskip
Let $M$ be a matroid over $V$ (described by its bases). The Tutte
polynomial is defined by
\[ t_M(x,y) = \sum_{X \subseteq V} (x-1)^{r(V)-r(X)}(y-1)^{n(X)} \]
where $n(X) = \min \{ |X \setminus Y| \mid Y \in M\}$ and $r(X)
= |X| - n(X)$ are the nullity and rank of
$X \subseteq V$ in $M$, respectively.

The following result is an extension of known results for
4-regular graphs \cite{Jaeger1990Transition} and for binary
matroids \cite[Proposition~4]{AignerPenroseBinMat} to matroids
in general.
\begin{Theorem} \label{thm:tutte_as_interlacep}
Let $M$ be a matroid over $V$. Then $Q_{[a,b,0]}(M)(y) =
a^{n(V)} b^{r(V)} t_M(1+\frac ab y, 1+\frac ba y)$.
\end{Theorem}
\begin{QProof}
We have $a^{n(V)} b^{r(V)} t_M(1+\frac ab y, 1+\frac ba y) =
a^{n(V)} b^{r(V)} \sum_{X \subseteq V} (\frac ab
y)^{r(V)-r(X)}(\frac ba y)^{n(X)}$ $= \sum_{X \subseteq V}$
$a^{|V\setminus X|} b^{|X|} y^{r(V)-r(X)+n(X)}$. It suffices to
show that $r(V)-r(X)+n(X) = d_M(X)$. We have $n(X) = \min \{ |X
\setminus Y| \mid Y \in M\}$. Hence $r(V)-r(X)+n(X) = r(V) -
|X| + 2n(X) = \min \{ (|Y| - |X| + |X \setminus Y|) + |X
\setminus Y| \mid Y \in M\}$ as $r(V)$ is the cardinality of the sets
$Y \in M$. Now $|Y| - |X| + |X \setminus Y| = |Y \setminus X|$
and thus we obtain $\min \{ (|Y| - |X| + |X \setminus Y|) + |X
\setminus Y| \mid Y \in M\} = \min \{ |Y \setminus X| + |X
\setminus Y| \mid Y \in M\} = \min \{ |X \sdif Y| \mid Y \in
M\} = d_M(X)$.
\end{QProof}

If we take $a=b=1$ in Theorem~\ref{thm:tutte_as_interlacep},
then we obtain the following corollary.
\begin{Corollary} \label{cor:tutte_as_interlacep}
Let $M$ be a matroid. Then $t_M(y,y) = q_1(M)(y-1)$.
\end{Corollary}

Hence, the Tutte polynomial on the diagonal is essentially $q_1(M)$ (in the sense that they both represent the same information) for the case where $M$ is a matroid described by its bases. The polynomial $R_M(x,y) = t_M(x+1,y+1)$ is known as the \emph{Whitney rank generating function}, as recalled in \cite[Section~15.4]{Welsh/MatroidBook}. Hence $R_M(y,y) = q_1(M)(y)$.

\section{Recursive Relations and \dmatroids}\label{sec:ss_polynomials_recursion}
In this section we show that each of the polynomials considered in
Section~\ref{sec:def_sets_interlacep} fulfill a specific
recursive relation when restricting to (vf-safe) \dmatroids.

Section~\ref{sec:equivalence_vf-safe} shows that
\dmatroids may be lifted to 2-matroids, and vf-safe \dmatroids
to tight 3-matroids.
This allows us to transfer properties of polynomials for multimatroids to similar properties for \dmatroids.

\subsection{The transition polynomial}

We first show that the polynomials $Q(M)$ for \vfsafe \dmatroids $M$ and $Q(Z)$ for tight $3$-matroids $Z$ are equivalent.
\begin{Theorem} \label{thm:mm_dm_pols_equal}
Let $M$ be a \vfsafe \dmatroid over $V$, $(U, \Omega)$ be a $(|V|,3)$-carrier, $\pi: U \to V$ be a projection of $(U, \Omega)$, and $\tau = (T_1,T_2,T_3)$ be a transversal $3$-tuple of $\Omega$. Then $Q(M)(\vec{a},\vec{b},\vec{c},y) = Q(\mathcal{Z}_{M,\tau})(\vec{x},y)$, where, for $u \in U$, $x_u = a_{\pi(u)}$ if $u \in T_1$, $x_u = b_{\pi(u)}$ if $u \in T_2$, $x_u = c_{\pi(u)}$ if $u \in T_3$.
\end{Theorem}
\begin{QProof}
Let $Z = \mathcal{Z}_{M,\tau}$. By Lemma~\ref{lem:null=dist}, $n_Z(T) = d_{M*\pi(T \cap T_2)\dual \pi(T \cap T_3)}$. Thus we have
\begin{eqnarray*}
Q(Z)(\vec{x},y) &=& \sum_{T \in \mathcal{T}(\Omega)} x_T y^{n_Z(T)}  = \sum_{T \in \mathcal{T}(\Omega)} a_{\pi(T \cap T_1)} b_{\pi(T \cap T_2)} c_{\pi(T \cap T_3)} y^{n_Z(T)} \\
&=& \sum_{(A,B,C) \in \mathcal{P}_3(V)} a_A b_B c_C
\, y^{d_{M*B\dual C}} = Q(M)(\vec{a},\vec{b},\vec{c},y).
\end{eqnarray*}
\end{QProof}

The ``$2$-matroid version'' of Theorem~\ref{thm:mm_dm_pols_equal} is obtained similarly.
\begin{Corollary} \label{cor:twomm_dm_pols_equal}
Let $M$ be a \dmatroid over $V$, $(U, \Omega)$ be a $(|V|,2)$-carrier, $\pi: U \to V$ be a projection of $(U, \Omega)$, and $\tau = (T_1,T_2)$ be a transversal $2$-tuple of $\Omega$. Then $Q(M)(\vec{a},\vec{b},\vec{c},y) = Q({\mathcal Z}_{M,\tau})(\vec{x},y)$, where, for $u \in U$, $x_u = a_{\pi(u)}$ if $u \in T_1$, $x_u = b_{\pi(u)}$ if $u \in T_2$, and $c_{v} = 0$ for all $v \in V$.
\end{Corollary}
\begin{QProof}
Since $c_v=0$ for all $v\in V$, it is sufficient to restrict the summation in the definition of $Q(M)$ to the triples $(A,B,C)\in \mathcal{P}_3$ with $C=\emptyset$. Other than this, the proof of this result is identical to the proof of Theorem~\ref{thm:mm_dm_pols_equal}.
\end{QProof}

We are now ready to state recursive relations for $Q(M)$.
\begin{Theorem} \label{thm:QM_recursive_relation}
Let $M$ be a \vfsafe \dmatroid and let $v \in V$.

\noindent
{\rm(1)} If $v$ is strongly nonsingular in $M$,
then
\[Q(M) = a_v Q(M\vertexrem v) + b_v Q(M*v\vertexrem v) +
c_v Q(M\dual v\vertexrem v).\]

\smallskip\noindent
{\rm(2)} Assume $v$ is not strongly nonsingular in $M$, and
let $\{\; (z_1,M_1),(z_2,M_2),(z_3,M_3) \;\} = \{\;(a_v,M\vertexrem
v), (b_v,M*v\vertexrem v), (c_v,M\dual v\vertexrem v)\;\}$. If
$M_1$ is empty, then $M_2 = M_3$ is nonempty and
\[Q(M) = (z_2 + z_3 + z_1 y)Q(M_2).
\]
\end{Theorem}
\begin{QProof}
Let $Z = \mathcal{Z}_{M,\tau}$ for some transversal 3-tuple $\tau = (T_1,T_2,T_3)$ for $(|V|,3)$-carrier $(U,\Omega)$ and projection $\pi: U \to V$. By Theorem~\ref{thm:mm_dm_pols_equal}, $Q(M)(\vec{a},\vec{b},\vec{c},y) = Q(Z)(\vec{x},y)$ for some suitable $U$-indexed vector $\vec{x}$. Let $\omega = \pi^{-1}(v)$.

Assume first that $v$ is strongly nonsingular in $M$. By Lemma~\ref{lem:minors_dm_mm}, $\omega$ is nonsingular in $Z$. By Theorem~\ref{thm:mm_recursive}, $Q(Z)(\vec{x},y) = \sum_{u \in \omega} x_{u} Q(Z|u)(\vec{x}',y) = a_v Q(Z|t_1)(\vec{x}',y) + b_v Q(Z|t_2)(\vec{x}',y) + c_v Q(Z|t_3)(\vec{x}',y)$ where $t_i \in T_i \cap \omega$ for all $i \in \{1,2,3\}$ and $\vec{x}'$ is obtained from $\vec{x}$ by removing the entries indexed by $\omega$. By Lemma~\ref{lem:minors_dm_mm}, for all $i \in \{1,2,3\}$, $Q(Z|t_i)(\vec{x}',y) = Q(\mathcal{Z}_{N_i,\tau'})(\vec{x}',y)$ with $N_1 = M\vertexrem v$, $N_2 = M*v\vertexrem v$, $N_3 = M\dual v \vertexrem v$, and $\tau' = (T_1 \setminus \omega, T_2 \setminus \omega, T_3 \setminus \omega)$. By Theorem~\ref{thm:mm_dm_pols_equal}, $Q(\mathcal{Z}_{N_i,\tau'})(\vec{x}',y) = Q(N_i)(\vec{a}',\vec{b}',\vec{c}',y)$ where the vectors $\vec{a}',\vec{b}',\vec{c}'$ are obtained from $\vec{a},\vec{b},\vec{c}$ by removing the entries indexed by $v$.

Assume now that $v$ is not strongly nonsingular in $M$. Assume that $M_1$ is empty. Let $\omega = \{t_1, t_2, t_3\}$ where for all $i \in \{1,2,3\}$, $t_i \in T_i$. By Lemma~\ref{lem:minors_dm_mm}, $\omega$ is singular in $Z$. Let $u \in \omega$ be singular in $Z$. By Theorem~\ref{thm:mm_recursive}, $Q(Z)(\vec{x},y) = (x_u y + \sum_{z \in \omega\setminus\{u\}}x_z)Q(Z|w)(x',y)$ for all $w \in \omega$. Since $Z|m = Z|n$ for all $m,n \in \omega$, we have by Lemma~\ref{lem:minors_dm_mm}, $\mathcal{Z}_{M_2,\tau'} = \mathcal{Z}_{M_3,\tau'}$ and thus $M_2 = M_3$ by Theorem~\ref{thm:vfsafe_equiv_t3matroid}. Thus, for all $w \in \omega$, $Q(Z|w)(\vec{x}',y) = Q(\mathcal{Z}_{M_2,\tau'})(\vec{x}',y)$. Finally, by Theorem~\ref{thm:mm_dm_pols_equal}, $Q(\mathcal{Z}_{M_2,\tau'})(\vec{x}',y) = Q(M_2)(\vec{a}',\vec{b}',\vec{c}',y)$ where the vectors $\vec{a}',\vec{b}',\vec{c}'$ are obtained from $\vec{a},\vec{b},\vec{c}$ by removing the entries indexed by $v$.
\end{QProof}

Note that, just like for Theorem~\ref{thm:mm_recursive}, the recursive relations of Theorem~\ref{thm:QM_recursive_relation} characterize $Q(M)(\vec{a},\vec{b},\vec{c},y)$.

We now state the obvious $2$-matroid version of Theorem~\ref{thm:QM_recursive_relation}.
\begin{Corollary} \label{cor:QM_recursive_relation_c}
Assume that $\vec{c}$ is the zero vector in $Q(M)$. Let $M$ be a
\dmatroid and let $u \in V$.

\noindent
{\rm(1)}  If $u$ is nonsingular in $M$,
then
 \[Q(M) = a_u Q(M\vertexrem u) + b_u Q(M*u\vertexrem u) .\]

\noindent
{\rm(2)} Assume $u$ is singular in $M$.
If $M\vertexrem u$ is empty, then $M*u\vertexrem u$ is nonempty and
\[Q(M) = (b_u + a_u y)Q(M*u\vertexrem u).
\]
A similar statement holds if $M*u\vertexrem u$ is empty.
\end{Corollary}

\subsection{The polynomials $q_i(M)$}

In the rest of this section we turn to various specializations of $Q(M)$.

We immediately obtain the following recursive relation for $q_1(M)$.
\begin{Theorem} \label{thm:ss_recur_interlacep}
Let $M$ be a \dmatroid. If $u \in V$ is nonsingular in $M$, then
\[ q_1(M) =
q_1(M\vertexrem u) + q_1(M*u\vertexrem u).\]
If every $v \in V$ is singular in $M$,
then $q_1(M) = (y+1)^n$ with $n = |V|$.
\end{Theorem}
\begin{QProof}
By Corollary~\ref{cor:QM_recursive_relation_c}, it suffices to consider the case where every $v \in V$ is singular in $M$, i.e., $M$ contains only one subset. Repeat the recursion $q_1(M)
= (1+y) q_1(M')$ where nonempty $M'$ is either $M\vertexrem u$ or
$M*u\vertexrem u$, until we reach $M_\emptyset
=(\emptyset,\{\emptyset\})$ which has $q_1(M_\emptyset) =1$.
\end{QProof}

\newcommand{\systeem}[2]{$(\,#1,\;\{#2\}\;)$}
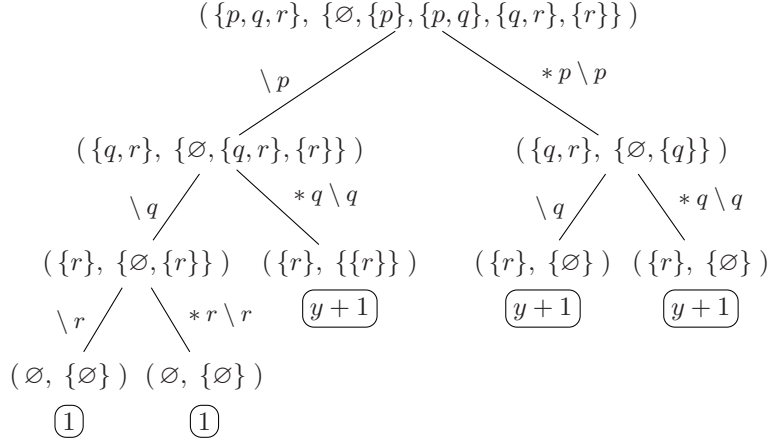
\begin{figure}[t]
\centerline{\unitlength 1.5mm
{\begin{picture}(70,41)(0,3)
\gasset{Nw=5,Nh=5,Nframe=n,Nfill=n,AHnb=0,ELpos=65}
  \node(+p-q)(46,20){ \systeem{\{r\}}{\emptyset} }
\node[Nframe=y,Nadjust=wh,Nmr=1](y)(46,16){$y+1$}
  \node(+p+q)(60,20){ \systeem{\{r\}}{\emptyset} }
\node[Nframe=y,Nadjust=wh,Nmr=1](y)(60,16){$y+1$}
  \node(-p-q-r)(04,10){ \systeem{\emptyset}{\emptyset} }
\node[Nframe=y,Nadjust=wh,Nmr=1](y)(04,6){$1$}
  \node(-p-q+r)(16,10){ \systeem{\emptyset}{\emptyset} }
\node[Nframe=y,Nadjust=wh,Nmr=1](y)(16,6){$1$}
  \node(-p-q)(10,20){ \systeem{\{r\}}{\emptyset,\{r\}} }
  \node(-p+q)(28,20){ \systeem{\{r\}}{\{r\}} }
\node[Nframe=y,Nadjust=wh,Nmr=1](y)(28,16){$y+1$}
  \node(root)(35,42){ \systeem{\{p,q,r\}} {\emptyset,\{p\},\{p,q\},\{q,r\},\{r\}} }
  \node(-p)(17,30){ \systeem{\{q,r\}}{\emptyset,\{q,r\},\{r\}}}
  \node(+p)(53,30){ \systeem{\{q,r\}}{\emptyset,\{q\}} }
\drawedge[ELside=r](-p-q,-p-q-r){\small${}\vertexrem r$}
\drawedge(-p-q,-p-q+r){\small${}*r\vertexrem r$}
\drawedge[ELside=r](-p,-p-q){\small${}\vertexrem q$}
\drawedge(+p,+p+q){\small${}*q\vertexrem q$}
\drawedge[ELside=r](+p,+p-q){\small${}\vertexrem q$}
\drawedge(-p,-p+q){\small${}*q\vertexrem q$}
\drawedge[ELside=r](root,-p){\small${}\vertexrem p$}
\drawedge(root,+p){\small${}*p\vertexrem p$}
\end{picture}}%
}
\caption{Recursive computation of $q_1(M)$.}\label{fig:recursive_q1}
\end{figure}

\begin{Example}
We recursively compute $q_1(M)$ using
Theorem~\ref{thm:ss_recur_interlacep} for \dmatroid $M =
\allowbreak (\{p,q,r\}, \allowbreak \{\emptyset, \{p\},
\{p,q\}, \{q,r\}, \{r\}\} )$ from
Example~\ref{ex:setsystem-orbit}. The computation tree is given
in Figure~\ref{fig:recursive_q1}. Recursion stops when the
\dmatroid is singular for all elements of the ground set. We verify that $q_1(M) = 3y+5$.
\end{Example}

Obviously $q_2(M)$ and $q_3(M)$ fulfill similar recursive
relations as $q_1(M)$ by applying the previous result to
$q_1(M\dual V)$ and $q_1(M+V)$, respectively. However, as the family of \dmatroids is not closed under loop complementation, we require that $M$ is vf-safe.
In this case we obtain, e.g., $q_2(M) =  q_2(M\dual u\vertexrem u) + q_2(M*u\vertexrem u)$.

We now consider the case where $M$ is normal, i.e., $\emptyset \in
M$. Note that $u \in V$ is nonsingular in a normal $M$ iff there is an $X \in
M$ with $u \in X$. Also, by Lemma~\ref{lem:ss_distance_prop}, $M$ is
normal iff $M+V$ is normal, and thus in this case $u$ is nonsingular in $M+V$ iff there is an $Y \in M+V$ with $u \in Y$.

We modify now Theorem~\ref{thm:ss_recur_interlacep}
such that each of the
``components'' of $M$ in the recursive equality are normal
whenever $M$ is normal. In this way we prepare for a
corresponding result on graphs.

\begin{Corollary} \label{cor:ss_recur_interlacep_series}
Let $M$ be a normal \vfsafe \dmatroid.\\
If $X \in M$ with $u \in X$, then
$$
q_1(M) = q_1(M\vertexrem u) + q_1(M*X\vertexrem u).
$$
If both $\{u,v\} \in M$ and $\{u\},\{v\} \not\in M$, then
\begin{eqnarray*}
q_2(M) &=& q_2(M*\{u,v\}\vertexrem \{u,v\}) +
q_2(M*\{u\}\dual \{v\}\vertexrem \{u,v\}) + \\
& & q_2(M\dual \{u\}\vertexrem u).
\end{eqnarray*}
Moreover, each of the given ``components'' is normal.
\end{Corollary}
\begin{QProof}
Since $q_1$ is invariant under pivot, $q_1(M*X\vertexrem u) =
q_1(M*\{u\}*(X\setminus\{u\})\vertexrem u) = q_1(M*\{u\}\vertexrem
u*(X\setminus\{u\})) = q_1(M*\{u\}\vertexrem u)$ and the result
follows by Theorem~\ref{thm:ss_recur_interlacep}.

Now consider $q_2$. We have the recursive formula
$q_2(M) = q_2(M*u\vertexrem u) + q_2(M\dual u\vertexrem u)$
(assuming $u$ is nonsingular in $M$).

One may easily verify that given
$\emptyset, \{u,v\} \in M$ and $\{u\},\{v\} \not\in M$, we have
$\emptyset, \{v\},\allowbreak \{u,v\} \in M\dual \{u\}$ and $\{u\} \not\in
M\dual \{u\}$. Therefore, $u$ is nonsingular in $M\dual \{u\}$.
Hence $q_2(M) = q_2(M*\{u\}\vertexrem u) +
q_2(M\dual \{u\}\vertexrem u)$. As we have $\emptyset, \{v\} \in
(M*\{u\}\vertexrem u) \dual \{v\}$, $v$ is nonsingular in this set system and therefore $q_2(M) = q_2(M*\{u,v\}\vertexrem \{u,v\}) +
q_2(M*\{u\}\dual \{v\}\vertexrem \{u,v\}) +
q_2(M\dual \{u\}\vertexrem u)$. We have that $\emptyset$ is
contained in each of the components $M*\{u,v\}\vertexrem \{u,v\}$,
$M*\{u\}\dual \{v\}\vertexrem \{u,v\}$, and
$M\dual \{u\}\vertexrem u$, hence they are all normal.
\end{QProof}

\subsection{The polynomial $Q_1(M)$}

In case $M$ is a normal \dmatroid, we find that $Q_1(M)$ may be
computed given $q_2(N)$ for all ``sub-set systems'' $N =
M\sub{X}$ of $M$.
\begin{Lemma} \label{Q1M_as_sum_q2M}
Let $M$ be a normal \dmatroid. Then $Q_1(M) = \sum_{X\subseteq V}
q_2(M\sub{X})$.
\end{Lemma}
\begin{QProof}
We start by observing that for any \dmatroid $M$ and $X\subseteq V$, $d_{M\sub{X}} = d_M$ provided that $M\sub{X}$ is nonempty. This follows from the fact that each minimal set (w.r.t.\ inclusion) of $M\sub{X}$ is a minimal set of $M$, and the observation from \cite[Property~4.1]{Bouchet_1991_67} that the minimal sets of \dmatroid $M$ are of equal cardinality.

Since $M$ is normal, $M\sub{X}$ is nonempty for all $X \subseteq V$.
We have therefore $q_2(M\sub{X}) = \sum_{Z \subseteq X} y^{d_{M\sub{X}+Z}(X)}$ for all $X \subseteq V$.
As $M$ is a \dmatroid, we have by the observation above, for all $Z \subseteq X$,
$d_{M\sub{X}+Z}(X) = d_{M+Z\sub{X}}(X) = d_{M+Z}(X)$. Hence $Q_1(M) = \sum_{X\subseteq V} q_2(M\sub{X})$.
\end{QProof}

Note that if $M$ is a \dmatroid (not necessarily normal), then
Lemma~\ref{Q1M_as_sum_q2M} can be applied, for any $Z\in M$, to the
normal \dmatroid $M*Z$ to obtain $Q_1(M*Z) = Q_1(M)$ (recall
that $Q_1$ is invariant under pivot).

From Theorem~\ref{thm:QM_recursive_relation}
we know that $Q_1(M)$ itself fulfills the
recursive relation
$$ Q_1(M) = Q_1(M\vertexrem u) + Q_1(M*u\vertexrem u) +
Q_1(M\dual u\vertexrem u),$$
where $M$ is vf-safe, and $u$ is strongly nonsingular in $M$.
%
If every $u \in V$ is not strongly nonsingular in $M$, then $Q_1(M) = (y+2)^n$ with $n=|V|$. Indeed, we have $Q_1(M) = (2+y) Q_1(M')$ where $M'$ is nonempty and equal to either $M\vertexrem u$ or $M*u\vertexrem u$. We have that every element of $M'$ is not strongly nonsingular in $M'$ (it is easy to see that if $v$ is strongly nonsingular in $M'$, then so is $v$ in $M$). By iteration, we reach $M_\emptyset =(\emptyset,\{\emptyset\})$ for which $Q_1(M_\emptyset) =1$.

The recursive relation given does not hold for
\dmatroids in general. Indeed, consider \dmatroid $M =
(V,2^V\setminus\{\emptyset\})$ with $V = \{1,2,3\}$. It is
shown in \cite{BH/NullityLoopcGraphsDM/10} that $M$ is not a
\vfsafe \dmatroid. We have $Q_1(M) = 13y+14$, while
$Q_1(M\vertexrem u) + Q_1(M*u\vertexrem u) + Q_1(M\dual
u\vertexrem u) = (3y+6) + (y+2)^2 + (y+2)^2 = 2y^2 + 11y + 14$.

\begin{figure}[t]
{\unitlength 1.5mm
\centerline{\begin{picture}(80,39)(-10,3)
\gasset{Nw=5,Nh=5,Nframe=n,Nfill=n,AHnb=0,ELpos=65}
  \node(-p-q)(00,10){$\{\emptyset, \{r\}\}$}
\node[Nframe=y,Nadjust=wh,Nmr=1](y)(00,06){$y+2$}
  \node(-p+q)(10,10){$\{\{r\}\}$}
\node[Nframe=y,Nadjust=wh,Nmr=1](y)(10,06){$y+2$}
  \node(-p*q)(20,10){$\{\emptyset\}$}
\node[Nframe=y,Nadjust=wh,Nmr=1](y)(20,06){$y+2$}
  \node(root)(30,40){$\{\emptyset, \{p\}, \{p,q\}, \{q,r\}, \{r\}\}$}
  \node(-p)(10,25){$\{\emptyset, \{q,r\}, \{r\}\}$}
  \node(+p)(30,25){$\{\emptyset, \{q\}\}$}
\node[Nframe=y,Nadjust=wh,Nmr=1](y)(30,21){$(y+2)^2$}
  \node(*p)(50,25){$\{\{q\}, \{q,r\}, \{r\}\}$}
  \node(*p-q)(40,10){$\{\{r\}\}$}
\node[Nframe=y,Nadjust=wh,Nmr=1](y)(40,06){$y+2$}
  \node(*p+q)(50,10){$\{\emptyset, \{r\}\}$}
\node[Nframe=y,Nadjust=wh,Nmr=1](y)(50,06){$y+2$}
  \node(*p*q)(60,10){$\{\emptyset\}$}
\node[Nframe=y,Nadjust=wh,Nmr=1](y)(60,06){$y+2$}
\drawedge[ELside=r](-p,-p-q){\small${}\vertexrem q$}
\drawedge[ELpos=70](-p,-p+q){\hspace{-2mm}\small${}*q\vertexrem q$}
\drawedge(-p,-p*q){\small${}\dual q\vertexrem q$}
\drawedge[ELside=r](root,-p){\small${}\vertexrem p$}
\drawedge(root,+p){\small${}*p\vertexrem p$}
\drawedge(root,*p){\small${}\dual p\vertexrem p$}
\drawedge[ELside=r](*p,*p-q){\small${}\vertexrem q$}
\drawedge[ELpos=70](*p,*p+q){\hspace{-2mm}\small${}*q\vertexrem q$}
\drawedge(*p,*p*q){\small${}\dual q\vertexrem q$}
\end{picture}}%
}
\caption{Recursive computation of $Q_1(M)$.}\label{fig:recursive_Q1}
\end{figure}
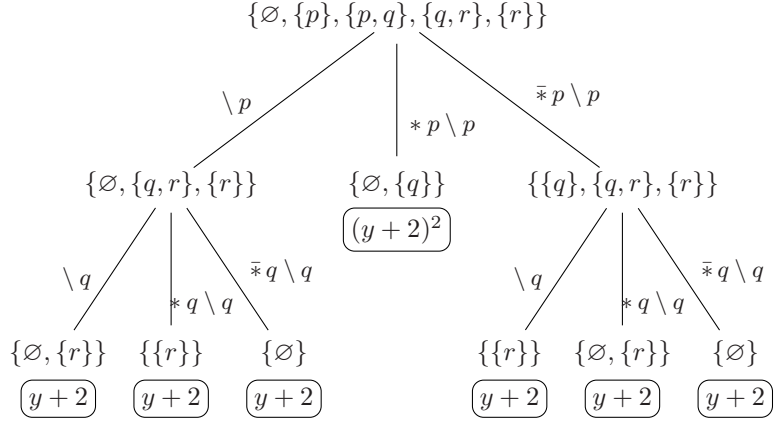

\begin{Example} \label{ex:recursive_Q1}
We recursively compute $Q_1(M)$ for \vfsafe \dmatroid
$M = (\{p,q,r\}, \{\emptyset, \{p\},\allowbreak \{p,q\}, \{q,r\},
\{r\}\} )$ from Example~\ref{ex:setsystem-orbit}. The
computation tree is given in Figure~\ref{fig:recursive_Q1},
where we have omitted the ground sets for visual clarity. The
recursion stops when every element of ground set of the \vfsafe \dmatroid under consideration is singular. We verify that $Q_1(M) = y^2 + 10y+16$.
\end{Example}

\subsection{The Tutte polynomial} \label{ssec:tutte_rec}
Let us reformulate Theorem~\ref{thm:ss_recur_interlacep} to
obtain a recursive characterization of $q_1(M)$ which can be
directly compared to that of the Tutte polynomial (for
matroids).

\begin{Corollary}
Let $M$ be a \dmatroid. Then
\begin{eqnarray*}
q_1(M) &=& \begin{cases}
(y+1)\, q_1(M \vertexrem u) & \mbox{if $M*u \vertexrem u$ is empty} \\
(y+1)\, q_1(M*u \vertexrem u) & \mbox{if $M \vertexrem u$ is empty} \\
q_1(M \vertexrem u) + q_1(M*u \vertexrem u) & {otherwise}
\end{cases}
\end{eqnarray*}
for all $u \in V$, and $q_1(M) = 1$ if $M =
(\emptyset,\{\emptyset\})$.
\end{Corollary}

We recall now the recursive relations of the Tutte polynomial
$t_M(x,y)$. Let $M$ be a matroid over $V$ (described by its
bases). Recall that $u$ is a coloop of $M$ iff $M \vertexrem u$
is empty, and $u$ is a loop iff $M*u\vertexrem u$ is
empty. Moreover, if $u \in V$ is not a coloop, then deletion of
matroids coincides with deletion of set systems, and if $u \in V$ is
not a loop, then contraction $M \slash u$ of $M$ by $u$
coincides with $M*u\vertexrem u$.
The notion of loop used here is not to be confused with loop
complementation and loops in graphs as we use in the rest of
the paper.

The Tutte polynomial $t_M(x,y)$ fulfills the following
characteristic relations:
\begin{eqnarray*}
t_M(x,y) &=& \begin{cases}
y\, t_{M \vertexrem u}(x,y) & \mbox{if $u$ is a loop} \\
x\, t_{M \slash u}(x,y) & \mbox{if $u$ is a coloop} \\
t_{M \vertexrem u}(x,y) + t_{M \slash u}(x,y) & \mbox{otherwise}
\end{cases}
\end{eqnarray*}
for all $u \in V$, and $t_M(x,y) = 1$ if $M =
(\emptyset,\{\emptyset\})$. Hence we find that the recursive
relations of $q_1(M)(y-1)$, restricted to matroids $M$, and
$t_M(x,y)$, for the case $x = y$, coincide.

Note that while matroids may be viewed as $1$-matroids, here we view the class of matroids also as a particular subclass of the class of even \dmatroids and therefore (essentially) as a particular subclass of the class of tight $2$-matroids. The class of vf-safe matroids may be viewed as a particular subclass of the class of tight $3$-matroids as well. Both the tight $2$-matroid viewpoint and the tight $3$-matroid viewpoint will be used in Section~\ref{sec:ss_poly_values} when we consider evaluations of the Tutte polynomial.

We remark that it is not possible to use the recursive relation of the Tutte polynomial to obtain a
full, two-variable generalization of the Tutte polynomial for
\dmatroids, as the outcome of the recursive relation would then depend
on the order of evaluations of the elements of $V$. For example, for the \dmatroid
$(\{a,b\},$ $\{\emptyset,$ $\{a\},\{a,b\}\}\,)$, the outcome of
the recursion might then be either $x+2$ or $y+2$. Of course,
we may use the extended recursive relation of
Corollary~\ref{cor:QM_recursive_relation_c} to cover a bigger
portion (instead of only the diagonal) of the two-variable
Tutte polynomial (for matroids), cf.\
Theorem~\ref{thm:tutte_as_interlacep}.

\subsection{Two-variable polynomials}

Using the generic recursive relation for \dmatroids,
Corollary~\ref{cor:QM_recursive_relation_c}, we formulate the
special case for the two-variable polynomial $\bar q(M)(x,y)$
defined in Section~\ref{ssec:two-var_&_tutte}.

\begin{Theorem} \label{thm:recur_two-var_interlace}
Let $M$ be a \dmatroid, and $u \in V$. If $u$ is nonsingular in $M$, then
\[ \bar q(M) =
\bar q(M\vertexrem u) + x \, \bar q(M*u\vertexrem u).\]

If $u$ is singular in $M$, then either
$M \vertexrem u$ is empty, and
$\bar q(M) = (x+y) \, \bar q(M *u \vertexrem u)$, or
$M*u \vertexrem u$ is empty, and
$\bar q(M) = (1+xy) \, \bar q(M\vertexrem u)$.
\end{Theorem}

If all $v \in V$ are singular in $M$, then $M$ contains a single element $X$,
and we can iterate the recursion to obtain
$\bar q(M) = (x+y)^{k}(1+xy)^{n-k}$ with $n = |V|$ and $k=|X|$.

Using Theorem~\ref{thm:recur_two-var_interlace} we derive a recursive relation in the style of the two-variable interlace polynomial for graphs \cite{ArratiaBS04}, see Section~\ref{sec:graphs}.

\begin{Lemma} \label{lem:umlooped-recursion}
Let $M$ be a \dmatroid, with $\emptyset, \{u,v\} \in M$,
while $\{u\},\{v\} \notin M$.
\\
Then $\bar q(M) =
\bar q(M\vertexrem u) + \bar q(M*\{u,v\}\vertexrem u) +
(x^2-1) \, \bar q(M*\{u,v\}\vertexrem \{u,v\} )
$
\end{Lemma}
\begin{QProof}
Obviously, $u$ is nonsingular in $M$. Hence by Theorem~\ref{thm:recur_two-var_interlace}, $\bar q(M) =
\bar q(M\vertexrem u) + x \, \bar q(M*u\vertexrem u)$.

(1) First assume that $v$ is nonsingular in $M*u\vertexrem u$. We have by Theorem~\ref{thm:recur_two-var_interlace}, $\bar q(M\vertexrem u) + x \, \bar q(M*u\vertexrem u) = \bar q(M\vertexrem u) + x \, \bar q(M*u\vertexrem \{u,v\}) + x^2 \, \bar q(M*\{u,v\}\vertexrem \{u,v\} )$.

Since $v$ is nonsingular in $M*u\vertexrem u$, $v$ is nonsingular in  $M*\{u,v\}\vertexrem u$. By Theorem~\ref{thm:recur_two-var_interlace} we find that
$
\bar q(M*\{u,v\}\vertexrem u)
=
\bar q(M*\{u,v\}\vertexrem u\vertexrem v) + x\, \bar q(M*\{u,v\}\vertexrem u *v\vertexrem v)
=
\bar q(M*\{u,v\}\vertexrem \{u,v\}) + x\, \bar q(M*u\vertexrem \{u,v\})
$. Thus $x\, \bar q(M*u\vertexrem \{u,v\}) = \bar q(M*\{u,v\}\vertexrem u) - \bar q(M*\{u,v\}\vertexrem \{u,v\})$, and we obtain the required equality.

(2) Now assume that $v$ is singular in $M*u\vertexrem u$. Note that $\emptyset \in M*\{u,v\} \vertexrem \{u,v\}$, which means that $(M*u\vertexrem u) *v \vertexrem v = M*\{u,v\} \vertexrem \{u,v\}$ is nonempty. We have by Theorem~\ref{thm:recur_two-var_interlace}, $\bar q(M\vertexrem u) + x \, \bar q(M*u\vertexrem u) = \bar q(M\vertexrem u) + x (x+y) \, \bar q(M*\{u,v\}\vertexrem \{u,v\})$.

Since $v$ is singular in $M*u\vertexrem u$, $v$ is singular in $M*\{u,v\}\vertexrem u$. By Theorem~\ref{thm:recur_two-var_interlace} we find that
$
\bar q(M*\{u,v\}\vertexrem u) =
(1+xy) \,\bar q(M*\{u,v\}\vertexrem \{u,v\})
$. Thus $xy \, \bar q(M*\{u,v\}\vertexrem \{u,v\}) = \bar q(M*\{u,v\}\vertexrem u) - \bar q(M*\{u,v\}\vertexrem \{u,v\})$, and we obtain the required equality.
\end{QProof}

\section{Polynomial Evaluations} \label{sec:ss_poly_values}
We consider now the polynomials of Section~\ref{sec:def_sets_interlacep} at particular values of $y$. We show in Section~\ref{sec:graph_poly} that various results in this section can be linked to (and are motivated by) known results corresponding to the case of graphs \cite{Aigner200411}. We also reconsider the evaluation of the Tutte polynomial $t_M$ at $(-1,-1)$ for quaternary matroids $M$ from \cite{DBLP:journals/jct/Vertigan98}.

\subsection{The polynomials $Q_1$ and $q_i$}
Note that we have $q_i(M)(1) = 2^n$ with $n = |V|$ and $i \in
\{1,2,3\}$. Moreover, $Q_1(M)(1) = \sum_{X,Y \subseteq V, X \cap Y =
\emptyset} 1 = 3^n$. Also, as $X \in M$ iff $d_M(X) = 0$, we have
that $q_1(M)(0)$ is equal to the number of sets in $M$.

\begin{Theorem} \label{thm:pol_eval_dmatroid}
Let $M$ be a \dmatroid.
\begin{enumerate}
\item \label{item:q1M_vmin1}
If $M$ is even and $|V| > 0$, then $q_1(M)(-1) = 0$.
\smallskip
\item \label{item:Q1M_vmin2}
If $M$ is \vfsafe and $|V| > 0$, then $Q_1(M)(-2)
= 0$.
\smallskip
\item \label{item:q1M_vmin2}
If $M$ is \vfsafe, then $q_1(M)(-2)
= (-1)^{|V|} (-2)^{d_{M\dual V}}$.
\end{enumerate}
\end{Theorem}
\begin{QProof}
By Corollary~\ref{cor:twomm_dm_pols_equal}, $q_1(M) = Q_{[1,1,0]}(M) = Q_1(\mathcal{Z}_{M,\tau})$ for suitable transversal $2$-tuple $\tau$ of some $(|V|,2)$-carrier $(U,\Omega)$ and projection $\pi$. Since $M$ is even, $\mathcal{Z}_{M,\tau}$ is tight by Subsection~\ref{ssec:tight_2m_even_dm}. By Theorem~\ref{thm:mm_Q_1-k}, we obtain (\ref{item:q1M_vmin1}).

By Theorem~\ref{thm:mm_dm_pols_equal}, $Q_1(M) = Q_1(\mathcal{Z}_{M,\tau})$ for suitable transversal $3$-tuple $\tau$ of some $(|V|,3)$-carrier $(U,\Omega)$ and projection $\pi$. Since $M$ is vf-safe, $\mathcal{Z}_{M,\tau}$ is tight by Theorem~\ref{thm:vfsafe_dmatroid_tight3matroid}. By Theorem~\ref{thm:mm_Q_1-k}, we obtain (\ref{item:Q1M_vmin2}).

Similarly, by Theorem~\ref{thm:mm_dm_pols_equal}, $q_1(M) = Q_{[1,1,0]}(M) = Q_1(\mathcal{Z}_{M,\tau}[U \setminus T_3])$ for suitable transversal $3$-tuple $\tau = (T_1,T_2,T_3)$ of some $(|V|,3)$-carrier $(U,\Omega)$ and projection $\pi$. Since $M$ is vf-safe, $Z = \mathcal{Z}_{M,\tau}$ is a tight $3$-matroid by Theorem~\ref{thm:vfsafe_dmatroid_tight3matroid}. By Theorem~\ref{thm:mmpol_res_1-k}, we obtain $Q_1(Z[U \setminus T_3]) = (-1)^{|\Omega|}(-2)^{n_Z(T_3)} = (-1)^{|V|} (-2)^{d_{M\dual V}}$ since $n_Z(T_3) = d_{M\dual V}$ by Lemma~\ref{lem:null=dist}. Thus we obtain (\ref{item:q1M_vmin2}).
\end{QProof}

Observe that similar properties hold for the other two polynomials
$q_i$, in particular $q_2(M)(-2) = (-1)^{|V|}(-2)^{d_M}$ for \vfsafe \dmatroid $M$.

Theorem~\ref{thm:pol_eval_dmatroid}(\ref{item:q1M_vmin1}) does not hold for \dmatroids in general:
for the \dmatroid $M=(V,2^V)$ the polynomial $q_1(M)$ equals (the constant)
$2^{|V|}$. Note that for the \vfsafe \dmatroid $M =
(\emptyset,\{\emptyset\})$ we have $Q_1(M)(-2) = 1$. Also, note
that for the (not \vfsafe) \dmatroid $M$ considered
before Example~\ref{ex:recursive_Q1},
we have $Q_1(M)(-2) = -12$.

\subsection{The Tutte polynomial}

It is a direct consequence of the recursive relation of the Tutte
polynomial that $t_M(0,0) = 0$ for matroid $M$,
except for $M = (\,\emptyset, \{ \emptyset \}\,)$.
This property may
now be seen as a special case of
Theorem~\ref{thm:pol_eval_dmatroid}(\ref{item:q1M_vmin1}) (as
every matroid is an even \dmatroid).

By Corollary~\ref{cor:tutte_as_interlacep} we may apply
Theorem~\ref{thm:pol_eval_dmatroid}(\ref{item:q1M_vmin2}) to obtain evaluations for the
Tutte polynomial.

\begin{Corollary} \label{cor:evaluations_tutte}
Let $M$ be a \vfsafe matroid. Then
$t_M(-1,-1) = (-1)^{|V|} (-2)^{d_{M\dual V}}$.
\end{Corollary}
The evaluation $t_M(-1,-1)$ for quaternary matroids was studied by Vertigan \cite{DBLP:journals/jct/Vertigan98}.
Let $A$ be a representation over $\gffour$ of a quaternary matroid $M$, and
let $C$ be the nullspace of $A$. Then the space $C \cap
C^{\perp}$ is called the \emph{bicycle space} of $C$. It turns out
that \cite[Proposition~3.5]{DBLP:journals/jct/Vertigan98} the dimension $d$ of the bicycle space is independent of
the representation $A$ over $\gffour$ of $M$, and $d$ is denoted by
$\mathrm{bd}(M,4)$ (where $4$ represents the field size). We
have $\mathrm{bd}(M,4) = \mathrm{bd}(M,2)$ for binary matroids
$M$. For binary matroids $M$, not only the dimension of the bicycle space, but the bicycle space $C \cap C^{\perp}$ itself is independent of the representation over $\two$.

It is shown in \cite[Theorem~3.4]{DBLP:journals/jct/Vertigan98} (by extending a result of \cite{Rosenstiehl1978195} for binary matroids) that
for quaternary matroids $M$, $t_M(-1,-1) = (-1)^{|V|}
(-2)^{\mathrm{bd}(M,4)}$. Consequently, by
Corollary~\ref{cor:evaluations_tutte} and the fact that quaternary matroids are \vfsafe, $d_{M\dual V} =
\mathrm{bd}(M,4)$ when $M$ is quaternary.
Note that $d_{M\dual V}$ provides a
natural extension of $\mathrm{bd}(M,4)$
to \vfsafe matroids $M$ that are not quaternary.

It is moreover shown in
\cite[Theorem~8.3]{DBLP:journals/jct/Vertigan98} that the
evaluation $t_M(-1,-1)$ is in some sense ``well behaved'' for
the larger class $\mathcal{N}$ of matroids that have no minors
isomorphic to $U_{2,6}$, $U_{4,6}$, $P_6$, $F_7^-$, or
$(F_7^-)^*$. We now illustrate this observation by applying
Corollary~\ref{cor:evaluations_tutte} to a non-quaternary
matroid.
\begin{Example}
The non-quaternary Steiner-system matroid $S(5,6,12)$ (see,
e.g., \cite{Oxley/MatroidBook-2nd} for a description of this
matroid) is \vfsafe \cite{BH/NullityLoopcGraphsDM/10}. One may
verify that $S(5,6,12) \dual V = S(5,6,12)$. Moreover,
$S(5,6,12)$ has rank $6$. Hence, by
Corollary~\ref{cor:evaluations_tutte}, $t_{S(5,6,12)}(-1,-1) =
2^{6}$.
\end{Example}

\subsection{Binary delta-matroids}
We now turn to binary \dmatroids.

For a $V\times V$-matrix $A$ (the columns and rows of $A$ are indexed by finite set $V$) over some field and $X \subseteq V$, $A[X]$ denotes the principal submatrix of $A$
with respect to $X$, i.e., the $X \times X$-matrix obtained
from $A$ by restricting to rows and columns in $X$. For a $V
\times V$-matrix skew-symmetric matrix $A$ over some field (we
allow nonzero diagonal entries over fields of characteristic $2$), the set system
$\mathcal{M}_A = (V,D_A)$ with $D_A = \{ X \subseteq V \mid
\det A\sub{X} \neq 0\}$ is a \dmatroid (see
\cite{bouchet1987}).
By convention, the nullity of the empty matrix is $0$, i.e.,
$n(A\sub{\emptyset}) = 0$.

A \dmatroid $M$ is called \emph{binary} if
$M = \mathcal{M}_A*X$ for some symmetric matrix $A$ over $\two$
and $X \subseteq V$. It is shown in \cite{bouchet1987} that a
matroid $M$ is binary in the above \dmatroid sense iff $M$ is
binary in the usual matroid sense. It is shown in
\cite{BH/NullityLoopcGraphsDM/10} that the family of binary
\dmatroids are closed under $*$ and $+$, and therefore every
binary \dmatroid is \vfsafe. In particular, every binary
matroid is \vfsafe.

The following is a reformulation of Corollary~4.2 of
\cite{TutteMartinOrientingVectors/Bouchet91} (and inspired by
Theorem~6 of \cite{Aigner200411}). The proof of
this reformulation is an appendix of Section~\ref{sec:graph_poly} as not all required notions are recalled at this point.

\begin{Proposition}[Corollary~4.2 of \cite{TutteMartinOrientingVectors/Bouchet91}] \label{prop:q1_2_Aigner}
Let $M$ be a binary \dmatroid. Then $q_1(M)(2) = k q_1(M)(-2)$ for some
odd integer $k$.
\end{Proposition}

The following example illustrates that Proposition~\ref{prop:q1_2_Aigner}
does not hold for \vfsafe \dmatroids in general.
\begin{Example} \label{ex:divide_bicycle}
Consider the $5$-point line $U_{2,5}$, which is a nonbinary \vfsafe matroid. We have $q_1(U_{2,5})(y) = y^3+6y^2+15y+10$. Hence, $q_1(U_{2,5})(2) = 9 \cdot 2^3$, while $q_1(U_{2,5})(-2) = -2^2$.
\end{Example}

\section{Pivot and Loop Complementation on Graphs}
 \label{sec:graphs}

In order to interpret the above results for graphs we recall in this
section the necessary notions and results from the
literature.
We recall that each graph has its characteristic \dmatroid,
defined via its adjacency matrix, and we recall two operations on graphs
that precisely match the \dmatroid operations pivot ${}*X$
and loop complementation ${}+X$.

We consider undirected graphs without parallel edges, but we do
allow loops.
For a graph $G=(V,E)$ we use $V(G)$ and $E(G)$ to denote its
set of vertices $V$ and set of edges $E$, respectively, and for
$x \in V$, $\{x\} \in E$ iff $x$ has a loop.

With a graph $G$ one associates its adjacency matrix $A(G)$,
which is a $V \times V$-matrix $\left(a_{u,v}\right)$ over
$\two$ with $a_{u,v} = 1$ iff $\{u,v\} \in E$ (we have $a_{u,u}
= 1$ iff $\{u\} \in E$). In this way, the family of graphs with
vertex set $V$ corresponds precisely to the family of
symmetric $V \times V$-matrices over $\two$. Therefore we
often make no distinction between a graph and its matrix, so,
e.g., by the nullity of graph $G$, denoted $n(G)$, we mean the
nullity $n(A(G))$ of its adjacency matrix (computed over
$\two$).
Also, we (may) write, e.g., $G\sub{X} = A(G)\sub{X}$ (for $X \subseteq V$), the subgraph of $G$ induced by $X$.
The graph $G \setminor X = G\sub{V\setminus X}$ is obtained from $G$ by deleting the vertices in $X$ with their incident edges. In case $X =
\{u\}$ is a singleton, we also write $G \vertexrem u$ to denote
$G \vertexrem \{u\}$.

Recall that we have defined the \dmatroid
$\mathcal{M}_A$ for a skew symmetric matrix $A$ as
$(V,D_A)$ with $D_A = \{ X \subseteq V \mid \det A\sub{X} \neq
0\}$,
so we set accordingly $\mathcal{M}_G = \mathcal{M}_{A(G)}$.
The construction adheres a correspondence between distance and nullity:
$d_{\mathcal{M}_G}(X) = n(G\sub{X})$ for $X \subseteq V$ (see
\cite{BH/NullityLoopcGraphsDM/10}).
Given $\mathcal{M}_G$ for some graph $G$, one can (re)construct
the graph $G$, see \cite[Property~3.1]{Bouchet_1991_67}.
In this way, the family of
graphs (with set $V$ of vertices) can be considered as a (strict) subset
of the family of binary \dmatroids (over set $V$).

For a graph $G$ and a set $X \subseteq V$, the graph obtained
after \emph{loop complementation} for $X$ on $G$, denoted by
$G+X$, is obtained from $G$ by adding loops to vertices $v \in
X$ when $v$ does not have a loop in $G$, and by removing loops
from vertices $v \in X$ when $v$ has a loop in $G$. Hence, if
one considers a graph as a matrix, then $G+X$ is obtained from
$G$ by adding the $V \times V$-matrix with elements $x_{i,j}$
such that $x_{i,i} = 1$ if $i \in X$ and $0$ otherwise. Note
that $(G+X)+Y = G+(X \sdif Y)$. It has been shown in
\cite{BH/PivotLoopCompl/09} that $\mathcal{M}_{G+X} =
\mathcal{M}_{G}+X$ for $X \subseteq V$ .

For a $V\times V$-matrix $A$ (over a field) and
$X \subseteq V$ such that $A\sub{X}$ is nonsingular, i.e.,
$\det A\sub{X} \neq 0$, the \emph{principal pivot
transform} (\emph{pivot} for short) of $A$ on $X$, denoted by
$A*X$, is defined as follows \cite{tucker1960}. Let $A = \begin{pmatrix}
P & Q \\
R & S
\end{pmatrix}$ with $P = A\sub{X}$. Then
$
A*X = \begin{pmatrix}
P^{-1} & -P^{-1} Q \\
R P^{-1} & S - R P^{-1} Q
\end{pmatrix}
$.
If $A$ is nonsingular, then $A*V = A^{-1}$.
The pivot operation is an
involution (operation of order $2$), and more generally, if
$(A*X)*Y$ is defined, then $A*(X \sdif Y)$ is defined and they
are equal.
If $A$ is skew-symmetric, then so is $A*X$. Thus,
this operation is defined on a graph $G$ via its adjacency matrix $A(G)$,
and yields a graph $G*X$.
Moreover, pivot operations on graphs and on their \dmatroids coincide:
 $\mathcal{M}_{G*X} = \mathcal{M}_{G}*X$
for $X \subseteq V$ if the left-hand side is defined \cite{bouchet1987} .

The pivots $G*X$ where $X$ is a minimal set of $\mathcal{M}_G
\backslash \{\emptyset\}$ (the set system obtained from
$\mathcal{M}_G$ by removing $\emptyset$) with respect to
inclusion are called \emph{elementary}.

It is noted in \cite{Geelen97} that any pivot on a graph can be
decomposed into a sequence of so-called elementary pivots,
which are operations ${}*X$ either on
a loop, $X = \{u\} \in E(G)$, or to an edge, $X =
\{u,v\} \in E(G)$, where (distinct) vertices $u$ and $v$ are
both non-loops.
The operation are known as \emph{local complementation}
and \emph{edge complementation}, respectively.
We do not recall (or use) their explicit
graph-theoretical definitions in this paper.
It can be found in, e.g., \cite{BH/PivotLoopCompl/09}.
Similar as for pivot we may define, e.g., ``dual local
complementation'' $G\dual \{u\} = G+\{u\}*\{u\}+\{u\}$ which is
identical to ``regular'' local complementation, except that it
is defined for a non-loop $\{u\}$ (instead of a loop).

For convenience, the next proposition summarizes the key known
results of this section, which we will use
frequently in the remaining part of this paper.

\begin{Proposition} \label{prop:known_results_graph_DM}
Let $G$ be a graph and $X \subseteq V$. Then the normal binary
\dmatroid $\mathcal{M}_G$ uniquely determines $G$ (and the
other way around). Moreover, $\mathcal{M}_{G*X} =
\mathcal{M}_{G}*X$ (if the left-hand side is defined),
$\mathcal{M}_{G+X} = \mathcal{M}_{G}+X$, and
$d_{\mathcal{M}_G}(X) = n(G\sub{X})$.
\end{Proposition}

\section{Graph Polynomials} \label{sec:graph_poly}
We now turn to graphs and reinterpret the above results on transition
polynomials for \dmatroids in this domain. We consider now
$q_i(M)$ for $i \in \{1,2,3\}$, $Q_1(M)$, $Q(M)$,
and (two-variable) $\bar q(M)$ for the
case $M = \mathcal{M}_G$ for some graph $G$. For notational
convenience we denote them by $q_i(G)$, $Q_1(G)$, $Q(G)$,
and $\bar q(G)$.

We obtain by Section~\ref{sec:def_sets_interlacep} and by
Proposition~\ref{prop:known_results_graph_DM} the following graph
polynomials.
\begin{eqnarray} \label{eqn:graph_interlacep_series}
q_1(G) &=& \sum_{X\subseteq V} y^{n(G\sub{X})} = \sum_{X\subseteq V} y^{n(G\vertexrem X)}, \\
q_2(G) &=& \sum_{X\subseteq V} y^{n(G+X)}, \\
q_3(G) &=& \sum_{X\subseteq V} y^{n(G+V\sub{X})}
\end{eqnarray}

Polynomial $q_1(G)(y-1)$ is the single-variable \emph{interlace
polynomial} $q(G)(y)$, which is defined in
\cite{Arratia2004199}. More generally, the two-variable
interlace polynomial \cite{ArratiaBS04} $q(G)(x,y) =
\sum_{X\subseteq V} (x-1)^{r(G\sub{X})} (y-1)^{n(G\sub{X})} =
\sum_{X\subseteq V} (x-1)^{|X|}
\left(\frac{y-1}{x-1}\right)^{n(G\sub{X})}$ is equal to $\bar
q(G)(x-1, \frac{y-1}{x-1})$. Moreover, the two-variable
polynomial
$\bar q(M*V \dual V) = \sum_{X \subseteq V} x^{|X|}
y^{d_{M*(V\setminus X)\dual X}} = \sum_{X \subseteq V} x^{|X|}
y^{d_{M+X}(V)}$ for $M=\mathcal{M}_G$
is equal to the
\emph{bracket polynomial} \cite{Traldi/Bracket1/09}, defined by
$b_{x}(G)(y) = \sum_{X\subseteq V} x^{|X|} y^{n(G+X)}$.
\footnote{Actually, the bracket polynomial contains another
variable $z$. However, as pointed out in
\cite{Traldi/Bracket1/09} we may assume without loss of
information either $x = 1$ or $z = 1$ (not both) --- for
convenience we choose $z = 1$.}

By the above, the definitions of $q_1(M)(y)$ and $q_2(M)(y)$ for
arbitrary set systems $M$ can be seen as generalizations of the
single-variable interlace polynomial and the single-variable
(i.e., the case $x=1$) bracket polynomial. From this point of
view we notice a close similarity between the single-variable
interlace polynomial and the single-variable bracket polynomial, as shown
in the next result (which follows directly from the definitions
of $q_2(M)$ and $q_3(M)$).

\begin{Theorem}
Let $G$ be a graph. We have (1) $q_1(G\dual V) = q_2(G)$, $q_1(G+V)
= q_3(G)$, $q_2(G*V) = q_3(G)$, assuming $G*V$ and $G\dual V$ are
defined, and (2) for $Y \subseteq V$, $q_1(G*Y) = q_1(G)$, $q_2(G+Y)
= q_2(G)$, $q_3(G\dual Y) = q_3(G)$, assuming $G*Y$ and $G\dual Y$
are defined.
\end{Theorem}

We now formulate the recursive relations from
Section~\ref{sec:ss_polynomials_recursion} concerning the $q_i(M)$
for normal set systems $M$ for the case that $M$ represents a graph
$G$.
\begin{Theorem} \label{thm:graph_recur_interlacep_series}
Let $G$ be a graph.
\begin{itemize}
\item Let $X \subseteq V$ such that $u \in X$ and $G*X$ is defined,
then
$$
q_1(G) = q_1(G\vertexrem u) + q_1(G*X\vertexrem u),
$$
\item
If $\{u,v\}$ is an edge of $G$ where both $u$ and $v$ do not
have loops, then
\begin{eqnarray*}
q_2(G)
&=& q_2(G*\{u,v\}\vertexrem \{u,v\}) +
q_2(G\dual \{v\}*\{u\}\vertexrem \{u,v\}) + \\
& & q_2(G\dual \{u\}\vertexrem u) \\
&=& q_2(G*\{u,v\}\vertexrem \{u,v\}) +
q_2(G*\{u,v\}\dual \{v\}\vertexrem \{u,v\}) + \\
& & q_2(G\dual \{u\}\vertexrem u), \hspace{2cm}\text{and}
\\[0.5em]
\bar q(G)
&=& \bar q(G \vertexrem u) +\bar q(G*\{u,v\}\vertexrem u ) +
    (x-1)^2\,\bar q(G*\{u,v\}\vertexrem \{u,v\}) \\
\end{eqnarray*}
\end{itemize}
\end{Theorem}
\begin{QProof}
The recursive relation for $q_1(G)$
follows directly
from Corollary~\ref{cor:ss_recur_interlacep_series}.
Likewise the relation for $\bar q(G)$ follows Lemma~\ref{lem:umlooped-recursion}.

We now consider $q_2(G)$. It is easy to see that $G*\{u,v\}$,
$G\dual \{v\}*\{u\}$, and $G\dual \{u\}$ are defined. Thus the
equality $q_2(G) = q_2(G*\{u,v\}\vertexrem \{u,v\}) +
q_2(G\dual \{v\}*\{u\}\vertexrem \{u,v\}) +
q_2(G\dual \{u\}\vertexrem u)$ follows now from
Corollary~\ref{cor:ss_recur_interlacep_series}. Moreover we have
$\mathcal{M}_G*\{u,v\}\dual \{v\}\vertexrem \{u,v\} =
\mathcal{M}_G*\{u\}*\{v\}\dual \{v\}\vertexrem \{u,v\} =
\mathcal{M}_G*\{u\}\dual \{v\}+\{v\}\vertexrem \{u,v\} =
\mathcal{M}_G*\{u\}\dual \{v\}\vertexrem \{u,v\} =
\mathcal{M}_G\dual \{v\}*\{u\}\vertexrem \{u,v\}$. As it is easy to
see that $G*\{u,v\}\dual \{v\}$ is defined, and so we have
$G*\{u,v\}\dual \{v\}\vertexrem \{u,v\} =
G\dual \{v\}*\{u\}\vertexrem \{u,v\}$. Consequently, we obtain the
other equality $q_2(G) = q_2(G*\{u,v\}\vertexrem \{u,v\}) +
q_2(G*\{u,v\}\dual \{v\}\vertexrem \{u,v\}) +
q_2(G\dual \{u\}\vertexrem u)$.
\end{QProof}

The recursive relation for the single-variable interlace
polynomial $q(G)$ in \cite[Section~4]{ArratiaBS04} is the
special case of the recursive relation for $q_1(G)$ in
Theorem~\ref{thm:graph_recur_interlacep_series} where the pivot
$G*X$ is elementary. Note that the recursive relation for
$q_1(M)$ where $M$ is a \dmatroid is a generalization of this
result. Also, it is shown in \cite{BH/PivotNullityInvar/09}
that the recursive relation for $q_1(G)$ in
Theorem~\ref{thm:graph_recur_interlacep_series} can be
generalized for arbitrary $V \times V$-matrices $A$. In this
way we find that $q_1(M)$ for \dmatroids and $q_1(A)$ for
matrices are ``incomparable'' generalizations (in the sense
that one is not more general than the other) of the graph
polynomial $q_1(G)$.

Of the two recursive relations for $q_2$ given in
Theorem~\ref{thm:graph_recur_interlacep_series}, the former seems
novel, while the latter is from
\cite[Theorem~1\emph{(ii)}]{Traldi/Bracket1/09}. The recursive relation for $\bar{q}$ given in
Theorem~\ref{thm:graph_recur_interlacep_series} corresponds to the recursive relation of the two-variable interlace polynomial \cite[Theorem~6]{ArratiaBS04}.

\begin{Remark}
The \emph{marked-graph bracket polynomial}
\cite{Traldi/Bracket2/09} for a graph $G$ is defined as
$\mathrm{mgb}_B(G,C) \allowbreak = \sum_{X\subseteq V} B^{|X|}
y^{n(G+X\sub{X \cup C})}$ with $C \subseteq V$ (the elements of
$C$ are the vertices of $G$ that are \emph{not} marked). For
the case $B = 1$, we have that $\mathrm{mgb}_B(G,C) = q_2(G*(V
\setminus C)) = q_3(G*C)$. Indeed, $q_3(M*C) = \sum_{X
\subseteq V} y^{d_{M*C\dual X}}$, and $d_{M*C\dual X} = d_{M*(C
\setminus X)*(C \cap X)\dual X} = d_{M*(C \setminus X)\dual
X+(C \cap X)} = d_{M*(C \setminus X)\dual X} = d_{M*(C
\setminus X)+X*X}$
$ = d_{M+X*(C \setminus X)*X} = d_{M+X*(X \cup
C)}\allowbreak = d_{M+X}(X \cup C)$ and so
$\mathrm{mgb}_1(G,C)$ is $q_3(M*C)$ where set system $M$
represents graph $G$ (i.e., $M = \mathcal{M}_G$). Therefore
$\mathrm{mgb}_1(G,C)$ can be seen as a ``hybrid'' polynomial
``between'' $q_2$ and $q_3$. The recursive relations for
$\mathrm{mgb}_B(G,C)$ deduced in \cite{Traldi/Bracket2/09} are
straightforwardly deduced from the recursive relations for
$q_i(G)$. Of course, one may also consider the hybrid
polynomials $q_1(G\dual C) = q_2(G\dual (V \setminus C))$ and
$q_1(G+C) = q_3(G+V \setminus C)$ for $C \subseteq V$ between
$q_1(G)$ and $q_2(G)$, and between $q_1(G)$ and $q_3(G)$,
respectively.
\end{Remark}

We now consider polynomial $Q_1(G)$. We have
\begin{eqnarray} \label{eqn:interlacep_Q}
Q_1(G) = \sum_{X\subseteq V} \sum_{Z\subseteq X} y^{n(G+Z\sub{X})}
\end{eqnarray}

By Lemma~\ref{Q1M_as_sum_q2M} we have $Q_1(G) = \sum_{X\subseteq V}
q_2(G\sub{X})$. By Theorem~\ref{thm:QM_recursive_relation} and the fact that
$Q_1(M)$ is invariant under pivot and loop complementation, we have
the following result.
\begin{Theorem} \label{thm:invar_pivot_loop_interlacep_graphs}
Let $G$ be a graph, and $Y \subseteq V$. We have $Q_1(G) = Q_1(G*Y)$
when $G*Y$ is defined, and $Q_1(G) = Q_1(G+Y)$. If $\{u,v\}$ is an
edge in $G$ with $u \not= v$ where both $u$ and $v$ are non-loop
vertices, then
\begin{eqnarray} \label{egn:recursive_rel_Q_loop}
Q_1(G) = Q_1(G\vertexrem u) + Q_1(G\dual \{u\}\vertexrem u) +
Q_1(G*\{u,v\}\vertexrem u).
\end{eqnarray}
\end{Theorem}
Note that the recursive relation of
Theorem~\ref{thm:invar_pivot_loop_interlacep_graphs} may be
easily modified for the cases where either $u$ or $v$ (or both) have
a loop. For example, if $u$ has a loop and $v$ does not have a loop
in $G$, then, since $Q_1(G) = Q_1(G+\{u\})$, we have $Q_1(G) =
Q_1(G\vertexrem u) + Q_1(G*\{u\}\vertexrem u) +
Q_1(G+\{u\}*\{u,v\}\vertexrem u)$ (as $G+\{u\}\dual \{u\}\vertexrem
u = G*\{u\}\vertexrem u$).

As $Q_1(G)$ is invariant under loop complementation, we may, as a
special case, consider this polynomial for \emph{simple graphs} $F$.
In this way, we obtain the recursive relation
\begin{eqnarray} \label{egn:recursive_rel_Q}
Q_1(F) = Q_1(F\vertexrem u) + Q_1(\mathrm{loc}_u(F)\vertexrem u) +
Q_1(F*\{u,v\}\vertexrem u),
\end{eqnarray}
where $\mathrm{loc}_u(F)$ is the operation that complements the
neighborhood of $u$ in $F$ without introducing loops (as $F$ is
simple, no loops are removed).

Polynomial $Q'(F) = \sum_{X\subseteq V} \sum_{Z\subseteq X}
(y-2)^{n(F+Z\sub{X})}$ has been considered in
\cite{Aigner200411} for simple graphs $F$\footnote{To ease
comparison, the formulation of the definition of $Q'(F)$ is
considerably modified with respect to the original
formulation in \cite{Aigner200411}.}. As the
difference in variable $y:=y-2$ in $Q_1(G)$ with respect to
$Q'(F)$ is irrelevant, this polynomial is essentially $Q_1(G)$
restricted to simple graphs. The
Equality~(\ref{egn:recursive_rel_Q}) is shown in
\cite{Aigner200411} from a matrix point-of-view. Similarly, the
single-variable case (case $u=1$) of the multivariate interlace
polynomial $C(F)$ of
\cite{DBLP:journals/combinatorics/Courcelle08} is (essentially)
$Q(G)$ restricted to simple graphs.

Finally, we may state the graph analogs of the results of
Section~\ref{sec:ss_poly_values}.
\begin{Theorem} \label{thm:summ_graph_results}
Let $G$ be a graph, and $n = |V|$. Then
\begin{eqnarray}
Q_1(G)(-2) &=& 0 \label{eqn:Q1G_vmin2} \quad \mbox{if $n>0$}\\
q_1(G)(-2) &=& (-1)^n (-2)^{n(G+V)} \label{eqn:q1G_vmin2} \\
q_2(G)(-2) &=& (-1)^n \\
q_1(G)(-1) &=& 0 \quad \mbox{if $n>0$ and $G$ has no loops} \label{eqn:q1G_vmin1} \\
q_1(G)(2) &=& k q_1(G)(-2) \quad \mbox{for some odd integer $k$}  \label{eqn:q1G_vplus2}
\end{eqnarray}
\end{Theorem}
\begin{QProof}
Equality~(\ref{eqn:Q1G_vmin2}) follows directly from
Theorem~\ref{thm:pol_eval_dmatroid}(\ref{item:Q1M_vmin2}).
As $d_{M\dual V} = d_{M+V*V+V} =
d_{M+V*V} = d_{M+V}(V)$, Equality~(\ref{eqn:q1G_vmin2}) follows
from Theorem~\ref{thm:pol_eval_dmatroid}(\ref{item:q1M_vmin2})
(and Proposition~\ref{prop:known_results_graph_DM}).
By the paragraph below
Theorem~\ref{thm:pol_eval_dmatroid} we have $q_2(G)(-2) = q_2(\mathcal{M}_G)(-2) = (-1)^n (-2)^{d_{\mathcal{M}_G}} = (-1)^n$.
Equality~(\ref{eqn:q1G_vplus2}) for $q_1(G)(2)$ follows from Proposition~\ref{prop:q1_2_Aigner}.
Finally, we have that $\mathcal{M}_G$ is an even \dmatroid iff
$G$ has no loops.
Hence Equality~(\ref{eqn:q1G_vmin1}) follows by
Theorem~\ref{thm:pol_eval_dmatroid}(\ref{item:q1M_vmin1}).
\end{QProof}

Equality~(\ref{eqn:Q1G_vmin2}) may also be
deduced from \cite[Section~2]{TutteMartinOrientingVectors/Bouchet91} through
the fundamental graph representation of isotropic systems
(and using that $Q_1(G)$ is invariant under loop complementation) ---
the result in \cite{TutteMartinOrientingVectors/Bouchet91} is in turn a
generalization of a result on the Martin polynomial for $4$-regular graphs.

Equality~(\ref{eqn:q1G_vmin2}) is proven in
\cite[Theorem~2]{Aigner200411},
\cite[Theorem~1]{DBLP:journals/ejc/BalisterBCP02}, and in
\cite[Corollary~1]{DBLP:journals/dm/Bouchet05}
for the case where $G$ does not have loops.
In fact in the first two of these proofs one
may recognize a type of ``$m,m,m+1$'' nullity result for matrices
which in the present paper is captured by the tightness property.

Equality~(\ref{eqn:q1G_vmin1}) is mentioned in \cite[Remark
before Lemma~2]{Aigner200411}. Note that
Equality~(\ref{eqn:q1G_vmin1}) does \emph{not} in general hold
for graphs with loops. Indeed, e.g., for the graph $G$ having
exactly one vertex $u$, where $u$ is looped, we have
$q_1(G)(-1) = 2$.


Equality~(\ref{eqn:q1G_vplus2}) is proven in
\cite[Theorem~6]{Aigner200411}, and in
\cite[Corollary~1]{DBLP:journals/dm/Bouchet05}
for the case where $G$ does not have loops.
The proof of Theorem~6 of \cite{Aigner200411} does not seem to be
easily extended from simple graphs to graphs with loops.

We remark that it is shown in \cite[Theorem~2]{Aigner200411} that $Q_1(G)(2) =
k2^n$ where $k$ is the number of induced Eulerian subgraphs.
This result cannot be generalized to \vfsafe \dmatroids as for
\vfsafe \dmatroid $M = ( \{q,r,s\}, \{ \{q\}, \{r\}, \{s\},
\{q,r\},$ $\{q,s\},$ $\{r,s\} \} )$ we have $Q_1( M ) =
9\cdot(y+2)$. Consequently, $Q_1(M)(2) = 9\cdot 2^2$ is not of
the form $k\cdot 2^n$.

Finally, a consequence of the proof of
\cite[Theorem~4.4]{bouchet1987} is that if $M$ is a binary
matroid, then $M = \mathcal{M}_G*X$ for some bipartite graph
$G$ (hence $G$ has no loops) with color classes $X$ and
$V\setminus X$. Conversely, if $G$ is a bipartite graph with
color classes $X$ and $V\setminus X$, then $M_G*X$ is a binary
matroid (and also $M_G*(V\setminus X)$). Hence we obtain the
following consequence of
Corollary~\ref{cor:tutte_as_interlacep}, which is essentially
shown in \cite[Theorem~3]{Aigner200411}.

\begin{Corollary}[Theorem~3 of \cite{Aigner200411}]
Let $G$ be a bipartite graph with color classes $X$ and $V\setminus X$,
and let $M = \mathcal{M}_G*X$ be a corresponding binary
matroid. Then $q_1(G)(y-1) = t_{M}(y,y)$.
\end{Corollary}

\subsection*{Appendix: proof of Proposition~\ref{prop:q1_2_Aigner}}

With the definition of the (single-variable) interlace
polynomial for graphs in place, we prove now
Proposition~\ref{prop:q1_2_Aigner} using Corollary~4.2 of \cite{TutteMartinOrientingVectors/Bouchet91} and Section~6 of \cite{Aigner200411}. We assume in this appendix that the reader is familiar with isotropic systems.
\begin{Proof-app}
Let $M$ be an binary \dmatroid, and let $X \in M$. Then $M*X =
\mathcal{M}_G$ for some graph $G$, and we have $q_1(M)(y) =
q_1(M*X)(y) = q_1(\mathcal{M}_G)(y) = q_1(G)(y) = q(G)(y+1)$,
where $q(G)(y)$ is the single-variable interlace polynomial. In
Section~6 of \cite{Aigner200411}, it is shown that, for graphs
$G$ without loops, $q(G)(y) = m(\mathcal{I}_G,A,y)$ where $m$
is the Martin polynomial for isotropic systems, $\mathcal{I}_G$
is a particular isotropic system depending on $G$ and $A$ is a
particular complete vector (we refer to
\cite{TutteMartinOrientingVectors/Bouchet91} for definitions of
these notions). It is straightforward to verify that the
reasoning of Section~6 of \cite{Aigner200411} holds essentially
unchanged for graphs $G$ where loops are allowed. The only
difference is that we may have either $(C_v,C_{\bar{v}}) =
(0,1)$ or $(C_v,C_{\bar{v}}) = (1,1)$ in the Claim of Section~6
of \cite{Aigner200411} instead of having only the former
possibility (and similarly for $D$), but this difference has no
effect on the conclusion that $\langle C,D \rangle = 1+1 = 0$.
By Corollary~4.2 of
\cite{TutteMartinOrientingVectors/Bouchet91}, for any isotropic
system $S$ and complete vector $A$ over the ground set of $S$,
$m(S,A,3) = k \cdot m(S,A,-1)$ for some odd integer $k$. We
conclude that $q_1(M)(2) = q(G)(3) = k q(G)(-1) = kq_1(M)(-2)$.
\qed\end{Proof-app}

\subsection*{Acknowledgements}
We thank Lorenzo Traldi for comments and corrections on the
paper, and for pointing out relevant literature to us.
We are indebted to the referees for their comments.
In particular, the paper benefited from the suggestion to move all basic arguments from the domain of delta-matroids to the domain of multimatroids, as this unified various arguments of Section~\ref{sec:ss_poly_values}.
R.B. is a postdoctoral fellow of the Research Foundation -- Flanders (FWO).

\bibliography{../geneassembly}

\end{document}